\documentclass[10pt,leqno]{article}

\usepackage[francais]{babel}
\usepackage{amsmath,amssymb,amscd}
\usepackage{a4wide}

\usepackage[latin1]{inputenc}

\usepackage[all]{xy} %%%pour une sortie pdf, il faut enlever dvips

\makeatletter

\@addtoreset{equation}{subsection}         %Met equation=0 lorsqu'une section commence

\makeatother

\newenvironment{paragr}[1][]{\refstepcounter{subsection} \noindent \textbf{\thesubsection . \ #1}}{\medskip}

\newenvironment{theoreme}{ \medskip\refstepcounter{theo}  \noindent\textbf{Th\'eor\`eme \thetheo}. ---\em}{\em \medskip}
\newenvironment{proposition}{\medskip\refstepcounter{theo}   \noindent\textbf{Proposition \thetheo}. ---\em}{\em\medskip}
\newenvironment{corollaire}{\medskip\refstepcounter{theo}  \noindent\textbf{Corollaire \thetheo}. ---\em}{\em\medskip}

\newenvironment{lemme}{\medskip\refstepcounter{theo}   \noindent\textbf{Lemme \thetheo}. ---\em}{\em\medskip}
\newenvironment{conjecture}{ \medskip\refstepcounter{theo}  \noindent\textbf{Conjecture \thetheo}. ---\em}{\em \medskip}

\newenvironment{definition}{\medskip\refstepcounter{theo}  \noindent\textbf{D\'efinition \thetheo}. ---}{\medskip}

\newenvironment{preuve}[1][]{\noindent \textbf{Démonstration.} #1 --- }{\hfill
  \ensuremath{\square} \medskip}

\newenvironment{remarque}{\medskip\refstepcounter{theo}  \noindent\textbf{Remarque \thetheo}. ---}{\medskip}

\DeclareMathOperator{\vol}{vol}

\DeclareMathOperator{\Hom}{Hom}

\DeclareMathOperator{\stab}{stab}

\DeclareMathOperator{\trace}{trace}

\newcommand{\ZZ}{\mathbb{Z}}

\newcommand{\NN}{\mathbb{N}}
\newcommand{\RR}{\mathbb{R}}
\newcommand{\AAA}{\mathbb{A}}
\newcommand{\CC}{\mathbb{C}}

\newcommand{\QQ}{\mathbb{Q}}

\newcommand{\Qp}{\mathbb{Q}_p}

\newcommand{\oc}{\mathcal{O}}

\newcommand{\Sc}{\mathcal{S}}

\newcommand{\ec}{\mathcal{E}}

\newcommand{\lc}{\mathcal{L}}

\newcommand{\fc}{\mathcal{F}}
\newcommand{\pc}{\mathcal{P}}

\newcommand{\ggo}{\mathfrak{g}}

\newcommand{\of}{\mathfrak{o}}

\newcommand{\mgo}{\mathfrak{m}}
\newcommand{\ngo}{\mathfrak{n}}

\newcommand{\pgo}{\mathfrak{p}}

\newcommand{\ugo}{\mathfrak{u}}

\newcommand{\al}{\alpha}

\newcommand{\la}{\lambda}

\newcommand{\back}{\backslash}

\newcommand{\Cc}{C_c^\infty}

\newcommand{\bg}{\langle}
\newcommand{\bd}{\rangle}

\newcommand{\eps}{\varepsilon}

\renewcommand{\leq}{\leqslant}
\renewcommand{\geq}{\geqslant}

\title{Sur une variante des troncatures d'Arthur}
\author{Pierre-Henri Chaudouard}
\date{}

\begin{document}
\maketitle
\begin{abstract}
Nous montrons que, pour une large classe de fonctions test, les contributions unipotentes dans la formule des traces pour $GL(n)$ sur un corps de nombre peuvent être obtenues à l'aide de fonctions zêta et d'intégrales de séries d'Eisenstein. La principale innovation est l'introduction d'une nouvelle troncature empruntée à un travail de Schiffmann sur les fibrés de Higgs.
  \end{abstract}

\renewcommand{\abstractname}{Abstract}
\begin{abstract}
  We show that, for a large class of test functions, the unipotent contributions in the trace formula for $GL(n)$ over a number field, can be obtained from zeta functions and integrals of Eisenstein series. The main innovation is a new truncation borrowed from  a work of Schiffmann on Higgs bundles.
\end{abstract}
\tableofcontents

\section{Introduction}

\begin{paragr}
  Soit  $F$ un corps de nombres et $\AAA$ son anneau d'adèles.
\end{paragr}

\begin{paragr}
Soit $n\geq 1$. Soit $G$ le groupe $GL(n)$ sur le corps $F$ et $\ggo$ son algèbre de Lie. Soit $\Sc(\ggo(\AAA))$ l'espace des fonctions complexes de Bruhat-Schwartz. Soit $\Cc(\ggo(\AAA))\subset \Sc(\ggo(\AAA))$ le sous-espace des fonctions à support compact.
\end{paragr}

\begin{paragr} Soit $f\in \Sc(\ggo(\AAA))$. La fonction $k$ définie sur le quotient $G(F)\back G(\AAA)$ par la formule
$$k(g)=k_f(g)=\sum_{X\in \ggo(F)} f(g^{-1}Xg)$$
est un analogue infinitésimal de la restriction à la diagonale du noyau automorphe. Soit $\oc$ l'ensemble des orbites géométriques définies sur $F$ de la représentation adjointe de $G$  sur  $\ggo$. On a un développement
$$k(g)=\sum_{\of\in \oc} k_\of(g)$$
où $k_\of(g)=\sum_{X\in \of(F)} f(g^{-1}Xg)$. Pour tout $\of\in \oc$ soit  $X\in \of(F)$ et $G_\of$ est le stabilisateur de $X$ dans $G$. On a alors un développement formel, analogue au côté géométrique de la formule des traces, 
$$\int_{G(F)\back G(\AAA)}  k(g) \,dg = \sum_{\of \in \oc}     \vol(G_\of(F)\back G_\of(\AAA)) \cdot    \int_{G_{\of}(\AAA)\back G(\AAA)} f(g^{-1}Xg)\,dg.
$$
  Malheureusement, dans le membre de droite aucun terme ne converge. On améliore légèrement  les choses en restreignant l'intégrale  au quotient $G(F)\back G(\AAA)^1$ de volume fini. Mais même avec cette restriction il existe des orbites dont le volume associé au centralisateur est infini et l'intégrale orbitale diverge. En suivant des méthodes d'Arthur (cf. \cite{ar1}), il est possible de  modifier  judicieusement $k$ en une fonction $k^\flat$ dont l'intégrale converge  (cf. \cite{PH1}).  Dans \cite{cuft}, on introduit pour tout $\of\in \oc$ une fonction modifiée   $k_\of^\flat$ dont la définition  est rappelée en \eqref{eq:kflat}. On  a alors
$$k^\flat(g)=\sum_{\of\in \oc} k_\of^\flat(g)$$
et chaque terme est alors intégrable sur  $G(F)\back G(\AAA)^1$. On obtient alors un développement 
$$
\int_{ G(F)\back G(\AAA)^1 } k^\flat(g)\, dg= \sum_{\of \in \oc} J_\of(f)
$$
où
$$J_\of(f)=\int_{ G(F)\back G(\AAA)^1 } k^\flat_\of(g)\, dg.
$$
Ce développement (établi dans \cite{cuft}, cf. aussi \cite{FiLa}) est un raffinement de l'analogue du développement géométrique usuel de la formule des traces, les termes étant paramétrés par les classes de conjugaison et non par les classes de conjugaison semi-simples comme chez Arthur.
\end{paragr}

\begin{paragr}[Le problème considéré.] --- La question à laquelle on s'intéresse est la suivante : peut-on trouver des formules explicites pour les contributions $J_\of(f)$ ? Précisons un peu la question. Pour des orbites $\of$ semi-simples régulières, la contribution $J_\of(f)$ est formellement simple : elle est donnée sous forme d'intégrales orbitales (éventuellement pondérées). En revanche, leur calcul (disons pour des fonctions test simples) est certainement très difficile et nous n'avons rien à ajouter ici. Lorsque $F$ est un corps de fonctions, ces intégrales comptent essentiellement  le nombre de points sur les corps finis de jacobiennes compactifiées de courbes singulières (courbes spectrales) associées à $\of$.

À l'opposé, les contributions nilpotentes ne sont pas sous une forme aussi aboutie (pour une réponse partielle, cf. \cite{scuft} et \cite{cuft}, voir également \cite{Hoff}). Pourtant elles doivent être de nature beaucoup plus élémentaire. Ainsi, lorsque $F$ est un corps de fonctions, la courbe spectrale est simplement un épaississement de la courbe projective lisse associée à $F$.

On relie aisément ces contributions nilpotentes à leur pendant unipotent dans la formule des traces pour $GL(n)$.  Cette question de l'explicitation des termes nilpotents est donc intéressante  lorsqu'on veut extraire d'une  seule formule des traces des informations de nature automorphe. Une autre motivation, lorsque $F$ est un corps de fonctions  et c'était notre motivation initiale, est d'obtenir un comptage de points de l'espace de modules des fibrés de Hitchin, en lien avec une conjecture de Hausel-Rodriguez-Villegas (cf. \cite{HRV}). Dans \cite{scfh}, ce comptage est directement relié à l'évaluation des distributions nilpotentes $J_\of(f)$ pour une fonction test très simple. Les constructions d'Arthur s'interprètent alors en termes de stabilité de fibrés vectoriels ou de fibrés de Hitchin (cf. \cite{Laf} et \cite{scfh}). Dans \cite{Schif} (cf. aussi \cite{MozS}), Schiffmann obtient une formule, certes compliquée,  pour ce comptage : pour cela, il tronque le champ algébrique des fibrés ou des fibrés de Hitchin en contrôlant la pente maximale d'un sous-fibré. Nous reprenons  ce point de vue dans le cadre des adèles des corps de nombres et de la formule des traces. On est ainsi conduit à introduire une fonction $k^\sharp_\of$ (cf. \eqref{eq:ksharp}) dont l'intégrale
$$J_\of^\sharp(f,s)=\int_{G(F)\back G(\AAA)} k^\sharp_\of(g) |\det(g)|^{s}\,dg
$$
converge absolument pour $s\in \CC$ de partie réelle $\Re(s)>0$. De plus, on a 
$$
\lim_{s\to 0}  s J^\sharp_\of(s)= J_\of(f)
$$
où  la limite est prise sur les $s\in \CC$ tels que $\Re(s)>0$ (cf. théorème \ref{thm:holo}). On introduit également l'intégrale orbitale tronquée
$$
J_\of^{G,\leq 0}(f,s)=\int_{G_\of(F)\back G(\AAA)} E^G(g) f(g^{-1}X  g) |\det(g)|^s \, dg.
$$
La fonction $E^G$ est la nouvelle fonction de troncature : de manière imagée, si l'on regarde  $g$ comme un fibré vectoriel, on a $E^G(g)=1$ si et seulement si la pente de tout sous-fibré est négative.

 L'intégrale $J_\of^{G,\leq 0}(f,s)$ converge absolument pour $\Re(s)>0$ (cf. proposition \ref{prop:IOtronq}). Sous une condition d'indépendance du terme constant (cf. la définition \ref{def:inv} de fonction \og presque invariante\fg{}), on a la formule suivante (cf. proposition \ref{prop:explicit})
 \begin{equation}
   \label{eq:Jsh}
   J_{\of}^{\sharp}(f,s)=\sum_{(M,\of') }\eps_M^G \frac{r_M}{|W(M,\of')|} J^{M,\leq 0}_{\of'}(f_M,s)
 \end{equation}
 où la somme est prise, à une certaine équivalence près,  sur les couples $(M,\of')$ formés d'un sous-groupe de Levi standard $M$  et d'une orbite $\of'$ dans l'algèbre de Lie de $M$ de sorte que l'induite de $\of'$ soit $\of$. Ici $f_M$ est le terme constant de $f$ (cf. §\ref{S:pinv}) et  $J^{M,\leq 0}_{\of'}$ est l'intégrale orbitale tronquée relative au sous-groupe de Levi $M$ (cf. \ref{S:IOLevi}). Le signe $\eps_M^G$ est défini en \eqref{eq:signe}. Les autres facteurs $|W(M,\of')|$ et $r_M$ sont des entiers naturels définis en \eqref{eq:WMof} et  \eqref{eq:rM}. Cette formule ramène donc essentiellement le calcul de $J_{\of}^{\sharp}(f,s)$ à celui de l'intégrale $J_\of^{G,\leq 0}(f,s)$.

Supposons désormais que  $\of$ est une orbite nilpotente. Dans ce cas, si $f$ est à support compact assez petit on a (cf. théorème \ref{thm:J0-J})
$$J^{G,\leq 0}_{\of}(f,s)= Z^G_{\of}(f,s)\cdot \theta_{\of}^G(s).
$$
Ici $Z^G_\of(f,s)$ est une fonction zêta de nature élémentaire (cf. §\ref{S:zeta}) qui converge et qui est holomorphe pour $\Re(s)>0$ (cf. proposition \ref{prop:zeta}). Par élémentaire, on entend que si $f$ est décomposable alors $Z^G_\of(f,s)$ est donné par un produit eulérien dont les facteurs presque partout s'expriment de manière combinatoire en termes des facteurs locaux de la fonction zêta du corps $F$. Le facteur $\theta_{\of}^G(s)$, défini à la section \ref{sec:Eis}, est indépendant de $f$ : c'est l'intégrale tronquée par la fonction $E^G$ d'une série d'Eisenstein élémentaire (on n'utilise rien de profond ; en particulier on ne sort pas du domaine de convergence de la série). Il converge pour $\Re(s)>0$ (cf. proposition \ref{prop:theta}). On obtient un énoncé similaire pour chaque sous-groupe de Levi $M$ : les objets affublés d'un exposant $M$ sont définis de la même façon mais relativement à $M$.

Voici la formule finale qu'on obtient pour $J_\of(f)$.

\begin{theoreme}(cf. théorème \ref{thm:residu})
  Soit $f\in \Cc(\ggo(\AAA))$ presque invariante au sens de la définition \ref{def:inv}. Lorsque le support de $f$ est assez petit, on  a,  pour toute orbite nilpotente $\of$,
$$J_\of(f)=\lim_{s\to 0}  s  \sum_{(M,\of')}\eps_M^G \frac{r_M}{|W(M,\of')|}Z^M_{\of'}(f_M,s)\cdot \theta_{\of'}^M(s).
$$
où la somme sur $(M,\of')$ est comme en \eqref{eq:Jsh} et la limite est prise sur les $s\in \CC$ tels que $\Re(s)>0$.
\end{theoreme}

Commentons les hypothèses et la nature de ce théorème. L'hypothèse \og presque invariante\fg{} est vérifiée pour des fonctions test intéressantes (cf.  remarque \ref{rq:inv}). L'hypothèse sur le support n'est pas si contraignante qu'elle paraît à première vue. Partant d'une fonction $f\in  \Cc(\ggo(\AAA))$, on peut toujours restreindre par une homothétie en une place son support et appliquer le théorème. En principe, on peut récupérer la valeur de $J_\of$ pour la fonction de départ en utilisant une propriété d'homogénéité des intégrales orbitales $J_\of(f)$ (cette propriété résulte des méthodes d'Arthur développées dans \cite{ar_unipvar}). Dans le théorème ci-dessus, tous les termes sont plutôt élémentaires de nature et, bien qu'on ne le fasse pas ici,  le facteur $\theta_{\of'}^M(s)$ peut en principe être calculé (cf. remarque \ref{rq:calcul-theta}). Néanmoins, les termes qui apparaissent ont, en général, des pôles d'ordre élevé ce qui rend le calcul pratique de la limite malaisé.

Pour terminer l'article, on donne pour une certaine fonction test, une formule conjecturale pour $J_\of(f)$. C'est l'analogue pour les corps de nombres d'un raffinement d'une conjecture de Hausel-Rodriguez-Villegas et Mozgovoy.
\end{paragr}

\begin{paragr}[Généralisations possibles.] --- La méthode semble, à première vue, très liée à $GL(n)$. En fait, on aurait pu tronquer le quotient $G(F)\back G(\AAA)^1$ par la condition dont l'analogue géométrique est que la pente de tout sous-fibré est inférieure à une constante $T$ donnée. Il s'agit ensuite de faire tendre $T$ vers l'infini. Sous cette forme, il nous semble que les méthodes de l'article devraient se généraliser si ce n'est à tout groupe réductif du moins aux groupes classiques. Une autre complication est qu'il faut tenir compte en général de la différence entre conjugaison rationnelle et conjugaison géométrique.  
\end{paragr}

\begin{paragr}[Notations.] --- \label{S:notations} Outre les notations qu'on a déjà introduites, on va utiliser les notations qui sont devenues standard depuis les travaux d'Arthur. Pour la commodité du lecteur, on en rappelle quelques-unes. Soit $T_0$ le sous-tore maximal diagonal de $G$ et $B$ le sous-groupe de Borel des matrices triangulaires supérieures. On qualifie de standard, resp. semi-standard, un sous-groupe parabolique défini sur $F$ qui contient  $B$, resp. $T_0$.  Soit $W$ le groupe de Weyl de $(G,T_0)$ qu'on identifie au sous-groupe des matrices de permutation. Pour tout $P$ sous-groupe parabolique   semi-standard, soit $M_P N_P$ sa décomposition de Levi où $M_P$ contient $T_0$ et $N_P$ est le radical unipotent de $P$. Dans la suite, le groupe $M_P$ est simplement appelé un sous-groupe de Levi de $G$. Soit $\pc(M)$ l'ensemble des sous-groupes paraboliques semi-standard $P$ de $G$  tels que $M_P=M$. On note par la lettre gothique correspondante $\pgo$, $\mgo_P$ et $\ngo_P$ l'algèbre de Lie des groupes $P$, $M_P$ et $N_p$.  Soit $K\subset G(\AAA)$ le sous-groupe compact maximal standard. On dispose de l'application d'Harish-Chandra $H_P: G(\AAA) \to a_P$ où $a_P=\Hom_\ZZ(X^*(P),\RR)$.  On pose 
$$\eps_P^G=(-1)^{\dim_\RR(a_P)-\dim_\RR(a_G)}.$$
et
\begin{equation}
  \label{eq:signe}
  \eps_{M_P}^G=\eps_P^G.
\end{equation}

On dispose d'ensembles de racines simples $\Delta_P$, resp. de poids simples $\hat{\Delta}_P$, qu'on voit comme des éléments de  $a_P^*=X^*(P)\otimes_\ZZ\RR$ et qui définissent des cônes  ouverts dans $a_P$ dont on note $\tau_P$, resp. $\hat{\tau}_P$, la fonction caractéristique. Soit $\rho_P$ la demi-somme des racines de $T_0$ dans $N_P$. Soit  $\Delta_0=\Delta_B$ et $a_0=a_B$ qu'on identifie naturellement à $\RR^n$. On a une inclusion naturelle $a_P\subset a_0$. On appelle type de $P$ la collection ordonnée des dimensions du gradué associé au drapeau de $F^n$  dont $P$ est le stabilisateur.

Les notations pour $G$ valent aussi lorsqu'on remplace $G$ par un sous-groupe parabolique $Q$. Dans ce cas, on ajoute un exposant $Q$ pour les objets relatifs à $Q$.

On fixe une mesure de Haar sur $G(\AAA)$. Soit $G(\AAA)^1$ le noyau de la valeur absolue adélique du déterminant. On munit $G(\AAA)$ et $G(\AAA)^1$ de mesures de Haar de sorte que le quotient soit muni, via son identification à $\RR_+^\times$, de la mesure $dt/t$ où $dt$ est la mesure de Lebesgue sur $\RR$. On fixe aussi des mesures  de Haar sur $M_P(\AAA)$, $N_P(\AAA)$ et $K$. On suppose que les volumes de $N_P(F)\back N_P(\AAA)$  et $K$ sont égaux à $1$. On suppose également que ces choix sont compatibles à la décomposition d'Iwasawa $G(\AAA)=M_P(\AAA)N_P(\AAA)K$.

Soit $Z_P=Z_{M_P}$ le centre de $M_P$. Soit $A_P$  le sous-tore $\QQ$-déployé maximal de $Res_{F/\QQ}(Z_P)$. Soit $A_P^\infty$ composante neutre de $A_P(\RR)$. Par une condition relative aux blocs linéaires,  on définit $M_P(\AAA)^1$. On  a alors une décomposition $M_P(\AAA)=M_P(\AAA)^1 A_P^\infty$. Là encore, on fixe des mesures de Haar compatibles à cette décomposition.

Pour alléger, on pose 
$$[G]=G(F)\back G(\AAA).
$$
\end{paragr}

\begin{paragr}[Remerciements.] ---  C'est un plaisir pour moi de remercier la fondation Simons et les organisateurs pour l'invitation au Simons Symposium sur la formule des traces et  le séjour extrêmement agréable au \emph{Schloss Elmau}.
Je remercie les projets Ferplay ANR-13-BS01-0012 et   Vargen ANR-13-BS01-0001-01 de l'ANR dont je fais partie. Je remercie  plus particulièrement l'Institut Universitaire de France qui me  fournit d'excellentes conditions de travail.
  \end{paragr}

\section{Variations sur les  constructions d'Arthur}

\begin{paragr}
  Dans cette section, on énonce et démontre de légères variations de résultats d'Arthur. Tout tourne autour de l'analogue pour les corps de nombres  de l'existence et l'unicité  de la  filtration de Harder-Narasimhan. La présentation s'inspire en particulier de  \cite{Behrend}.
\end{paragr}

\begin{paragr}[Degré d'instabilité.]
 Soit  $g\in G(\AAA)$ et $Q$ un sous-groupe parabolique standard. On définit \emph{le degré d'instabilité} de $g$ suivant $Q$ par la formule suivante
$$
\deg^Q_i(g)= \max_{ (P,\delta)  }\bg \rho_P^Q, H_P(\delta g) \bd 
$$
où  $(P,\delta)$ parcourt les couples formé d'un sous-groupe parabolique standard $P\subset Q$ et d'un élément $\delta\in P(F)\back Q(F)$. La borne supérieure des $\bg \rho_P^Q, H_P(\delta g) \bd $ est finie et atteinte comme il résulte de la théorie de la réduction (cf. lemme  5.1 de \cite{ar1}). Ce degré est invariante par $Z_Q(\AAA)$.

\begin{lemme}\label{lem:une-equiv}
Les trois conditions suivantes sont équivalentes
\begin{enumerate}
\item $$
\deg^Q_i(g)\leq 0
$$
\item pour tout sous-groupe parabolique $P\subset Q$, tout $\delta \in P(F)\back Q(F)$ et tout $\varpi\in \hat{\Delta}_P^Q$, on a 
$$\bg \varpi,H_P(\delta g) \bd \leq 0.$$
\item pour tout sous-groupe parabolique $P\subset Q$ maximal, tout $\delta \in P(F)\back Q(F)$ et tout $\varpi\in \hat{\Delta}_P^Q$, on a 
$$\bg \varpi,H_P(\delta g) \bd \leq 0.$$
\end{enumerate}
\end{lemme}

\begin{preuve}
On a $1\Rightarrow 3$ car si $P\subset Q$ est maximal et si $\varpi$ est l'unique élément de $\hat{\Delta}_P^Q$, on a $\rho_P^Q=c\cdot\varpi$ pour une constante $c>0$.

On a  $3\Rightarrow 2$ car tout $\varpi\in \hat{\Delta}_P^Q$ détermine un sous-groupe parabolique $P\subset P'\subset Q$ maximal et 
$$\bg \varpi,H_P(\delta g) \bd =\bg \varpi,H_{P'}(\delta g) \bd $$
ne dépend que de la classe de $\delta$ dans $P'(F)\back G(F)$.

Enfin  $2\Rightarrow 1$ car le vecteur $\rho_P^Q$ est une combinaison linéaire à coefficients positifs d'éléments de $\hat{\Delta}_P^Q$. 
\end{preuve}

 Lorsque $Q=G$, on pose
$$\deg_i(g)=\deg^G_i(g)$$
et on parle simplement de degré d'instabilité. On dit que $g$ est \emph{semi-stable} si et seulement si 
$$
\deg_i(g)\leq 0.
$$

\end{paragr}

\begin{paragr}[Couple canonique.] ---

\begin{definition}\label{def:cano}
On dit qu'un couple   $(P,\delta)$ formé d'un sous-groupe parabolique standard $P$ et d'un élément $\delta\in P(F)\back G(F)$ est \emph{canonique} s'il vérifie les deux conditions suivantes
\begin{enumerate}
\item $$
\bg \rho_P^G, H_P(\delta g) \bd =\deg_i(g)
$$
\item pour tout sous-groupe parabolique $Q$ contenant $P$ tel que $ \bg \rho_Q^G, H_Q(\delta g) \bd =\deg_i(g)$ on a $Q=P$.
\end{enumerate}
\end{definition}

\begin{lemme}\label{lem:reformulation}
  Soit  $g\in G(\AAA)$. Alors $(P,\delta)$ est un  couple canonique  de $g$ si et seulement si les deux conditions suivantes sont vérifiées
  \begin{enumerate}
  \item on a 
$$
\deg^P_i(\delta g)\leq 0
$$
\item Pour tout $\al\in \Delta_P$, on a 
$$
\bg \al, H_P(\delta g) \bd >0
$$
  \end{enumerate}

\end{lemme}

\begin{preuve}
  Soit $(P,\delta)$ un couple canonique de $g$.  Montrons que les conditions 1 et 2 sont nécessaires. On a  pour tout sous-groupe parabolique $Q\subset P$ et tout $\eta\in Q(F)\back P(F)$
\begin{eqnarray}
  \label{eq:=1} \bg \rho_Q^P, H_Q(\eta \delta g) \bd &=&   \bg \rho_Q^G, H_Q(\eta \delta g) \bd- \bg \rho_P^G, H_P( \delta g) \bd  \\
\nonumber &=& \bg \rho_Q^G, H_Q(\eta \delta g) \bd-\deg_i(g)\leq  0 
\end{eqnarray}
d'après la condition 1 de la définition \ref{def:cano}. On obtient ainsi   la première condition d'après le lemme \ref{lem:une-equiv}. 

Soit $\al\in \Delta_P$ et $P\subsetneq R \subset G$ le sous-groupe parabolique minimal défini par la condition $\Delta_P^R=\{\al\}$. On a alors
\begin{eqnarray}
 \label{eq:=2} \bg \rho_P^R, H_P( \delta g) \bd  &=&   \bg \rho_P^G, H_Q( \delta g) \bd-  \bg \rho_R^G, H_R(\delta g) \bd  \\
\nonumber&=& \deg_i(g)- \bg \rho_R^G, H_R(\eta \delta g) \bd >  0 
\end{eqnarray}
L'inégalité est stricte par la condition 2 de la définition \ref{def:cano} ($P$ est maximal parmi les paraboliques qui vérifient  $\deg_i(g)=\bg \rho_R^G, H_P(\delta g) \bd$). On obtient donc la seconde condition  puisque $\rho_P^R$ et $\al$ sont égaux à un coefficient multiplicatif près, qui est strictement positif.

Inversement, il est aisé de partir des égalités \eqref{eq:=1} et  \eqref{eq:=2}   pour obtenir la réciproque.
\end{preuve}
\end{paragr}

\begin{paragr}[Existence et unicité.] --- C'est l'objet de la proposition suivante.

\begin{proposition}\label{prop:HN}
  Pour tout $g\in G(\AAA)$, il existe un unique couple canonique.
\end{proposition}

\begin{preuve}
  L'existence est claire. Il s'agit de prouver l'unicité. Soit $(P,\delta)$ et $(Q,\eta)$ deux couples canoniques. On a alors $\delta \eta^{-1}\in P(F)\back G(F) / Q(F)$. On peut donc supposer qu'on peut trouver des relèvements à $G(F)$, notés encore  $\delta$ et $\eta$, tels que $\delta\eta^{-1}=w\in W$. Quitte à changer $g$ en $\eta g$, on peut supposer qu'on a $\eta=1$ et $\delta=w$. Ainsi, on considère les couples canoniques $(P,w)$ et $(Q,1)$. Les conditions 1 et 2 du lemme \ref{lem:reformulation} (cf. aussi lemme \ref{lem:une-equiv}) impliquent les conditions suivantes.
  \begin{itemize}
  \item Pour tout sous-groupe parabolique standard $R\subset P$, pour tout $\delta\in R(F)\back P(F)$ et tout $ \varpi\in \hat{\Delta}^P_R$, on a $\bg \varpi, H_R(\delta wg) \bd \leq 0$ ;
  \item Pour tout $\al\in \Delta_P$, on a $
\bg \al, H_P(wg) \bd >0.
$
  \end{itemize}
Posons $P_1=w^{-1}Pw$. Les conditions ci-dessus impliquent
\begin{itemize}
  \item Pour tout sous-groupe parabolique semi-standard $R\subset P_1$  et tout $ \varpi\in \hat{\Delta}^{P_1}_R$, on a $\bg \varpi, H_R(g) \bd \leq 0$ ;
  \item Pour tout $\al\in \Delta_{P_1}$, on a $\bg \al, H_{P_1}(g) \bd >0.
$
  \end{itemize}
L'élément $g$ vérifie les conditions analogues où l'on remplace $(P_1,w)$ par (Q,1). 

 Pour tout sous-groupe parabolique semi-standard $P'$, soit $C_{P'}(g)\subset a_0$ la partie formée des $H\in a_0$ qui satisfont les deux conditions suivantes :
\begin{itemize}
  \item Pour tout sous-groupe parabolique semi-standard $R\subset P'$ et tout $ \varpi\in \hat{\Delta}^{P'}_R$, on a $\bg \varpi, H+H_{P'}(g) \bd \leq 0$ ;
  \item Pour tout $\al\in \Delta_{P'}$, on a $
\bg \al, H+H_{P'}(g) \bd >0.
$
  \end{itemize}
Lorsque $P'$ parcourt l'ensemble des sous-groupes paraboliques contenant $T_0$, la collection des $C_{P'}(g)$ forme une partition de $a_0$ (cf. la discussion de \cite{localtrace} p.22).  D'après ce qui précède, on a 
$$0\in C_{P_1}(g) \cap C_{ Q}(g)
$$
donc  $P_1= Q$. Donc $P$ et $Q$ sont standard et conjugués donc égaux et donc $w\in P(F)\cap W$ ce qu'il fallait voir.
\end{preuve}
\end{paragr}

\begin{paragr}[Identités d'Arthur.] --- \label{S:IdAr}Soit $F^P$ la fonction caractéristique des $g\in G(\AAA)$ tels que $\deg^P_i(g)\leq 0$.

  \begin{proposition}\label{prop:ecriture} Pour tout sous-groupe parabolique $Q$ de $G$, on a 
    \begin{enumerate}
    \item On a 
$$1=\sum_{P\subset Q} \sum_{\delta \in P(F)\back Q(F)}F^P(\delta g) \tau_P^Q(H_P(\delta g))
$$
\item 
$$
F^Q(g)=\sum_{P\subset Q} \eps_P^Q \sum_{\delta \in P(F)\back Q(F)} \hat{\tau}_P^Q(H_P(\delta g))
$$
    \end{enumerate}
  \end{proposition}

  \begin{preuve} Ces assertions se ramènent à des assertions analogues sur le facteur de Levi standard $M_Q$ de $Q$. Ce dernier est un produit de groupes généraux linéaires et les fonctions qui apparaissent sont des produits sur les facteurs linéaires. Par récurrence, on suppose donc les assertions connues pour les sous-groupes paraboliques propres de $G$.
    L'assertion 1 pour $Q=G$ n'est qu'une reformulation de l'existence et l'unicité d'un couple canonique.  D'après l'assertion 1 pour $Q=G$ et l'assertion 2 pour $Q\subsetneq G$, on a
    
    \begin{eqnarray*}
      F^G(g)&=&1 -\sum_{Q\subsetneq G} \sum_{\delta \in Q(F)\back G(F)} F^Q(\delta g) \tau_Q^G(H_Q(\delta g))\\
&=& 1 -\sum_{Q\subsetneq G}  \sum_{\delta \in Q(F)\back G(F)}   \tau_Q^G(H_Q(\delta g))  \sum_{P\subset Q} \eps_P^Q \sum_{\eta \in P(F)\back Q(F)} \hat{\tau}_P^Q(H_P(\delta \eta g))\\
&=&  1 -\sum_{P\subsetneq G}\sum_{\delta \in P(F)\back G(F)}    \sum_{P\subset Q\subsetneq G}  \eps_P^Q \hat{\tau}_P^Q(H_P(\delta  g))  \tau_Q^G(H_Q(\delta g))\\
&=&   1 +\sum_{P\subsetneq G}\sum_{\delta \in P(F)\back G(F)}      \eps_P^G \hat{\tau}_P^G(H_P(\delta  g)) \\
    \end{eqnarray*}
car, d'après le lemme de Langlands, on a
\begin{equation}
  \label{eq:langlands}
  \sum_{P\subset Q\subset G}  \eps_P^Q\hat{\tau}_P^Q(H_P(\delta  g))  \tau_Q^G(H_Q(\delta g))=0.
\end{equation}
pour $P\subsetneq G$ (une référence récente est  \cite{LabWal}  proposition 1.7.2).
  \end{preuve}

  \begin{remarque}
    Dans \cite{ar1}, Arthur introduit à l'aide de la théorie de la réduction, une fonction $F^Q(\cdot,T)$ qui dépend d'un paramètre $T$. Notre fonction $F^Q$ n'est pas la valeur en $T=0$ de la fonction d'Arthur. Cependant, pour $T$ \og positif\fg{}, on a 
$$F^Q \leq F^Q(T).$$
En particulier, $F^Q$ hérite de certaines propriétés de support de $F^Q(T)$.  Ainsi la fonction $F^G$, qui est évidemment $Z_G(\AAA)$-invariante, est la fonction caractéristique d'un compact de $Z_G(\AAA)G(F)\back G(\AAA)$.
  \end{remarque}
\end{paragr}

\section{Une autre troncature}

\begin{paragr}
Dans cette section, on introduit une troncature   analogue à celle qui consiste à limiter le champ des fibrés vectoriels sur une courbe en imposant que la pente  d'un sous-fibré soit toujours négative.
\end{paragr}

\begin{paragr}
Soit $\Pi\subset a_0$ l'ensemble des poids dominants des représentations standard de $G$ à savoir les puissances extérieures de la  représentation naturelle. On a donc
$$
\Pi^G=\Pi=\{(1,0,\ldots,0), \ldots, (1,1,\ldots,1)\}.
$$
On appelle sous-groupe parabolique extrémal de $G$ un sous-groupe parabolique standard qui est soit égal à $G$ soit maximal. Les sous-groupes  paraboliques extrêmaux sont en bijection avec les éléments de $\Pi^G$ si l'on associe au sous-groupe parabolique standard $P$ de type $(r,n-r)$ le poids
$$\begin{array}{cc}\varpi_r=&(\underbrace{1,\ldots,1},\underbrace{0,\ldots,0}) \\&  {r }\,\,\, \text{      }  \, \,  {n-r }
  \end{array}.
  $$
Pour tout sous-groupe parabolique standard $P$ de $G$ on pose alors
$$\varpi_P=\frac1{n_1} \varpi_{n_1}
$$
où $(n_1,n_2,\ldots,n_k)$ est le type de $P$. 
Si $P\subsetneq G$ est maximal de type $(r,n-r)$, on a 
$$\hat{\Delta}_P^G=\{r(\varpi_P-\varpi_G )  \}.
$$
Soit $\chi_P^G$ la fonction caractéristique des $H\in a_0$ tels que 
$$
\bg \varpi_P, H \bd \leq 0.
$$
Soit $E^G$ l'application à valeurs dans $\{0,1\}$ définie pour tout $g\in G(\AAA)$ par
\begin{equation}
  \label{eq:EQ}
  E^G(g)= \sum_{P\subset G} \sum_{\delta \in P(F)\back G(F)}F^P(\delta g) \tau_P^G(H_P(\delta g)) \chi_P^G(H_P(\delta g))
\end{equation}
où la somme est prise sur les sous-groupes paraboliques standard de $G$. D'après la proposition \ref{prop:ecriture}, il y a au plus un terme non nul dans la somme ci-dessus.

\begin{remarque}
Il est immédiat sur la définition que cette application se descend au quotient $G(F)\back G(\AAA)/K$. Contrairement à la fonction $F^G$ définie au §\ref{S:IdAr}, elle n'est pas $Z_G(\AAA)$-invariante. Sur le sous-groupe $G(\AAA)^1$, les fonctions $F^G$ et $E^G$ coïncident. Il n'en est pas de même en général.
\end{remarque}
\end{paragr}

\begin{paragr} La définition précédente se généralise de la manière suivante. Soit $Q\subset G$ un  sous-groupe parabolique standard. Soit 
  $(l_1,\ldots,l_k)$  son type. Pour tout tout sous-groupe parabolique standard $P\subset Q$,  soit $$(n_1,\ldots, n_{r_1},n_{{r_1}+1},\ldots, n_{r_2}, n_{r_2+1},\ldots,n_{r_k}).
$$
On suppose que pour $1\leq i \leq k$, on a l'égalité  $l_i=n_{r_{i-1}+1} + \ldots+  n_{r_i}$ où l'on pose $r_0=0$. Soit $\chi_P^Q$ la fonction caractéristique des $H\in a_0$ tels que 
$$
 \sum_{j=l_1+\ldots+l_{j}+1}^{l_1+\ldots +l_{j}+n_{r_j+1} }  \bg e_j^* ,H\bd \leq 0
$$
pour $0\leq j \leq k-1$.

Soit $E^Q$ l'application à valeurs dans $\{0,1\}$ définie pour tout $g\in G(\AAA)$ par
\begin{equation}
  \label{eq:EQQ}
  E^Q(g)= \sum_{P\subset Q} \sum_{\delta \in P(F)\back Q(F)}F^P(\delta g) \tau_P^Q(H_P(\delta g)) \chi_P^Q(H_P(\delta g)).
\end{equation}

On obtient ainsi une application sur le quotient $ M_Q(F)N_Q(\AAA)\back G(\AAA)/ K$. Si l'on écrit $g=mnk$ selon la décomposition d'Iwasawa $G(\AAA)=M_Q(\AAA)N_Q(\AAA)K$ et $m=(m_1,\ldots, m_k ) $ selon l'isomorphisme $M_Q=GL(l_1)\times \ldots \times GL(l_k)$ on a 
$$E^Q(g)=  \prod_{i=1}^k E^{GL(l_i)}(m_i).$$
On note $E^{M_Q}$ la fonction sur $M_Q(\AAA)$ donné par  le membre de droite ci-dessus.
\end{paragr}

\begin{paragr} Les sous-groupes extrêmaux de $G$ sont munis d'un ordre total : on a  $P \leq Q$ si leurs types respectifs notés $(n_1,n_2)$ et $(m_1,m_2)$ vérifient $n_1\leq m_1$.

  \begin{proposition}\label{prop:mumax}
    Soit $g\in G(\AAA)$. 
    \begin{enumerate}
        \item La borne supérieure, prise sur les couples $(P,\delta)$ formés d'un sous-groupe parabolique extrémal et  d'un élément $\delta\in P(F)\back G(F)$,
$$\sup_{(P,\delta)} \bg \varpi_P,H_P(\delta g)\bd
$$
est finie et atteinte. 
\item Soit $P$ le plus grand sous-groupe parabolique extrémal (pour la relation d'ordre définie ci-dessus) tel qu'il existe $\delta\in P(F)\back G(F)$ de sorte que $(P,\delta)$ soit un couple qui atteint cette borne.  
Alors le couple $(P,\delta)$ est unique. Plus précisément, soit $(Q,\gamma)$ le couple canonique de type noté $(n_1,\ldots,n_r)$. Alors $P$ est de type $(n_1, n_2+\ldots+n_r)$ et $\delta$ est la projection de $\gamma$ sur $P(F)\back G(F)$.
\end{enumerate}
\end{proposition}

\begin{preuve}
L'assertion 1 résulte de la théorie de la réduction (cf. lemme 5.1 de \cite{ar1}).

Prouvons l'assertion 2. Soit $(P,\delta)$ comme dans l'assertion 2. Soit $(a,b)$ le type de $P$. En utilisant l'existence d'un couple canonique pour le facteur $GL(b)$, on voit qu'il existe $(Q,\gamma)$ avec $Q\subset P$ de type $(a,b_1,\ldots,b_k)$ et $\gamma \in P(F)\delta$ qui vérifient les deux conditions suivantes

\begin{enumerate}
\item[(A)] pour tout sous-groupe parabolique maximal $R$ de $Q$ de type $(a,b_1',\ldots,b'_{k+1})$, tout $\eta\in R(F)\back Q(F)$  et tout $\varpi\in \hat{\Delta}_R^{Q}$
$$
\bg \varpi,H_P(\eta \gamma g)\bd \leq 0
$$
\item[(B)] pour tout $\al\in \Delta_Q^P$ on a $\bg \al, H_Q(\gamma g)\bd >0$.
\end{enumerate}
Soit $R$ le sous-groupe parabolique standard maximal de $Q$ de type $(a_1,a_2,b_1,\ldots,b_k)$ avec $a_1+a_2=a$. Soit $R'$ le sous-groupe parabolique maximal de $G$ de type $(a_1,n-a_1)$. Pour tout élément $\eta\in Q(F)$, on a 

$$\bg \varpi_{R'},H_{R'}(\eta \gamma g)\bd \leq \bg \varpi_{P},H_{P}(\gamma g)\bd
$$
Donc pour $\varpi$ l'unique élément de $\hat{\Delta}_{R}^Q$, on a 
\begin{equation}
  \label{eq:3eg}
   \bg  \varpi , H_{R'}(\eta \gamma g) \bd =  \bg   a_1(\varpi_{R'}-\varpi_{P} ) , H_{R'}(\eta \gamma g) \bd \leq 0
 \end{equation}
 
 Soit $S$ le sous-groupe parabolique standard de $G$ de type $(a+b_1,b_2,\ldots,b_k)$. Soit $S'$  le sous-groupe parabolique standard maximal de $G$ de type $(a+b_1,b_2+\ldots+b_k)$. On a donc
$$\bg \varpi_{S'},H_{S'}(\gamma g)\bd \leq \bg \varpi_{P},H_{P}(\gamma g)\bd
$$
c'est-à-dire
$$0\leq \bg \varpi_{P}- \varpi_{S'},H_{Q}(\gamma g)\bd= \frac{b_1}{a+b_1}  \bg \al ,H_{Q}(\gamma g)\bd
$$
où $\al$ est l'unique élément de $\Delta_Q^S$.
donc par minimalité de $P$, on a 
\begin{equation}
  \label{eq:4eg}
 \bg \al ,H_{Q}(\gamma g)\bd >0. 
\end{equation}
En combinant les propriétés (A) et (B), les inégalités \eqref{eq:3eg} et  \eqref{eq:4eg} avec le lemme \ref{lem:reformulation}, on voit que $(Q,\gamma)$ est le couple canonique de $g$.
\end{preuve}
\end{paragr}

\begin{paragr}[Quelques corollaires.] ---

\begin{corollaire}\label{cor:car}
  On a $E^G(g)=1$ si et seulement si pour tout couple $(P,\delta)$ formé d'un sous-groupe parabolique extrémal et  d'un élément $\delta\in P(F)\back G(F)$ on a 
$$ \bg \varpi_P,H_P(\delta g)\bd \leq 0.
$$
\end{corollaire}

\begin{preuve}
  C'est immédiat : d'une part on a défini $E^G$ en terme du couple canonique (cf. \eqref{eq:EQ}) et d'autre part la proposition \ref{prop:mumax} traduit l'inégalité cherchée en terme du couple canonique. 
\end{preuve}

\begin{corollaire}\label{cor:encad}
  Soit $P$ un sous-groupe parabolique standard de type noté $(n_1,\ldots, n_r)$. Soit $p\in P(\AAA)$. On écrit $p=mn$ avec $m=(m_1,\ldots,m_r)\in M_P(\AAA)\simeq GL(n_1,\AAA)\times\ldots\times GL(n_r,\AAA)$
$$E^P(p)\leq E^G(p)\leq E^{GL(n_1)}(m_1).
$$
\end{corollaire}

\begin{preuve}
  L'inégalité  $E^G(p)\leq E^{GL(n_1)}(m_1)$ est aisé à obtenir compte tenu du corollaire \ref{cor:car}. Pour obtenir l'inégalité $E^P(p)\leq E^G(p)$ on peut raisonner par récurrence sur le nombre de blocs de $M_P$. Le cas crucial est le cas où $P$ a deux blocs, c'est-à-dire $P$ est maximal, cas qu'on considère maintenant.  On suppose qu'on $E^P(p)=1$. On veut montrer que $E^G(p)=1$. Soit $Q$ un sous-groupe parabolique extrémal de $G$ et $\delta\in Q(F)\back G(F)$. D'après le corollaire \ref{cor:car}, il s'agit de voir qu'on a 
  \begin{equation}
    \label{eq:but}
     \bg \varpi_Q,H_Q(\delta p)\bd \leq 0.
   \end{equation}
   La cas $Q=G$ est évident. On suppose donc $Q\subsetneq G$ c'est-à-dire $Q$ est maximal. On peut remplacer $p$ par $\gamma p$ avec $\gamma\in P(F)$ sans changer la condition $E^P(p)=1$. Quitte à faire ce changement, on peut et on va supposer que $\delta$ est représenté par $w\in W$ et que le sous-groupe parabolique $M_P\cap (w^{-1}Qw)$ est standard. On sait que
$$H_{w^{-1}Bw}(p)-H_B(p)$$
est une combinaison linéaire à coefficients positifs d'éléments $\al\in \Delta_B$ tels que $w\al$ est négatif pour $B$ (cf. \cite{dis_series} lemme 3.6).
On a donc
$$H_B(wp)=w\cdot H_{w^{-1}Bw}(p)= w\cdot H_B(p)+\sum_{\al \in \Delta_B} n_\al \al$$
 où $n_\al\leq 0$. Il s'ensuit qu'on a
 \begin{equation}
   \label{eq:inegw}
   \bg \varpi_Q,H_Q(w p)\bd \leq \bg w^{-1}\cdot \varpi_Q,H_B(p)\bd =  \bg w^{-1}\cdot \varpi_Q,H_B(m)\bd   
 \end{equation}
 Introduisons alors le sous-groupe parabolique $ R= (M_P\cap (w^{-1}Qw)) N_P \subset P$. Écrivons $M_P=GL(a)\times GL(b)$ et $R\cap M_P=R_1\times R_2$. On a alors les vecteurs $\varpi_{R_1}$ et $\varpi_{R_2}$ relatifs à $GL(a)$ et $GL(b)$ qu'on voit naturellement comme des vecteurs de $a_0$. Le majorant de \eqref{eq:inegw} s'écrit alors $\bg \varpi_{R_1}+\varpi_{R_2}, H_B(m)\bd$. Pour $i=1,2$,  les $R_i$ sont extrémaux et  on a  donc  $\bg \varpi_{R_i}, H_B(m)\bd\leq 0$. On obtient bien \eqref{eq:but}.
\end{preuve}
\end{paragr}

\section{Noyaux modifiés}\label{sec:integ-tronquee}

\begin{paragr}
  Dans cette section, suivant \cite{cuft}, on introduit des \og noyaux \fg{}  $k_\of^\flat$ dont l'intégrale sur $G(F)\back G(\AAA)^1$ est la contribution de l'orbite $\of$ qui nous intéresse. En exploitant la troncature décrite dans la section précédente, on définit un nouveau  noyau  $k_\of^\sharp$ dont on relie l'intégrale à celle de $k_\of^\flat$.
\end{paragr}

\begin{paragr} Soit $f\in \Sc(\ggo(\AAA))$. Pour tout sous-groupe parabolique standard $P$ de $G$ et toute orbite $\of\in \oc$, soit
$$k_{P,\of}(f,g)=\sum_{X\in \mgo_P(F), I_P^G(X)=\of} \int_{\ngo_P(\AAA)} f(g^{-1}(X+U)g)\, dU,
$$
où l'on note $I_P^G(X)$ l'induite de la $M_P$-orbite de $X$ (cf. \cite{cuft} §2.9). La mesure sur $\ngo_P(\AAA)$ est normalisée de sorte que le réseau  $\ngo_P(F)$ soit de covolume $1$.
Soit
$$k_{P,\of}^{\leq 0}(g)= E^{P}(g) k_{P,\of}(f,g).
$$
\end{paragr}

\begin{paragr}
 Soit 
\begin{equation}
  \label{eq:kflat}
k^\flat_\of(f,g)=\sum_{P} \eps_P^G  \sum_{\delta \in P(F)\back G(F)} \hat{\tau}_P(H_P(\delta g)) k_{P,\of}(f,\delta g)
\end{equation}
et
\begin{equation}
  \label{eq:ksharp}
  k^{\sharp}_\of(f,g)=\sum_{P} \eps_P^G  \sum_{\delta \in P(F)\back G(F)}\hat{\tau}_P(H_P(\delta g)) k_{P,\of}^{\leq 0}(f,\delta g)
\end{equation}
les sommes portent sur les sous-groupes paraboliques standard $P$ de $G$.
\end{paragr}

\begin{paragr}[Convergence d'une intégrale.] --- Pour tous sous-groupes paraboliques standard $P_1\subset P_2$, soit
  \begin{equation}
    \label{eq:sigma12}
    \sigma_{P_1}^{P_2}(H)=\sum_{P_2\subset P} \eps_{P_2}^{P}\tau_{P_1}^{P}(H) \hat{\tau}_{P}(H) 
  \end{equation}
  et
$$ k_{P_1,\of}^{P_2}(f,g)=\sum_{P_1\subset P\subset P_2} \eps_P^G   k_{P,\of}(f,g)
$$
où $P$ est un sous-groupe parabolique. La fonction  $\sigma_{P_1}^{P_2}$ est à valeurs dans $\{0,1\}$. Pour une autre description de cette fonction, on renvoie à \cite{ar1} lemme 6.1.  

Le théorème \ref{thm:cv} ci-dessous est un cas particulier du corollaire 3.2.2 de \cite{cuft}. En fait, on peut en  donner une preuve plus directe fondée sur la proposition \ref{prop:cv} ci-dessous. Cette dernière  se démontre avec les méthodes de \cite{cuft} section 3 (cf. en particulier la proposition 3.6.1)  à l'aide de la proposition \ref{prop:ecriture}.

 \begin{proposition}\label{prop:cv} 
  Pour tous sous-groupes paraboliques standard $P_1\subset P_2$, l'intégrale

$$\int_{P_1(F)\back G(\AAA)^1}F^{ P_1}(g) \sigma_1^2(H_{P_1}( g))k_{P_1,\of}^{P_2}(f,g) \, dg
  $$
est absolument convergente.
  \end{proposition}

  \begin{theoreme}\label{thm:cv}
    \begin{enumerate}
    \item     L'intégrale 
\begin{equation}
  \label{eq:Jof}
  J_\of(f)=
\int_{G(F)\back G(\AAA)^1 } k^\flat_\of(f, g)\, dg
\end{equation}
est absolument convergente.
\item De plus, on a
$$J_\of(f)=\sum_{P_1\subset P_2}\int_{P_1(F)\back G(\AAA)^1}F^{ P_1}(g) \sigma_1^2(H_{P_1}( g))k_{P_1,\of}^{P_2}(f,g) \, dg
$$
où la somme porte sur les sous-groupes paraboliques standard $P_1$ et $P_2$ tels que $P_1\subset P_2$.
\end{enumerate}
\end{theoreme}

\begin{preuve}
Par des manipulations standard (qui repose sur la proposition \ref{prop:ecriture} et la combinatoire des   pp. 41-43 de \cite{Ar-intro}), l'expression $k^\flat_\of(f, g)$ est égale à 
$$\sum_{P_1\subset P_2} \sum_{\delta\in P_1(F)\back G(F)}  F^{ P_1}(\delta g) \sigma_1^2(H_{P_1}( \delta g)) k_{P_1,\of}^{P_2}(f,\delta g) .$$
Le théorème est donc une conséquence directe de la proposition \ref{prop:cv}.
  \end{preuve}
\end{paragr}

\begin{paragr}[Un énoncé auxiliaire.] --- Pour tous sous-groupes paraboliques standard $P_1\subset P_2$, soit 
$$ k_{P_1,\of}^{P_2,\leq 0}(g)=\sum_{P_1\subset P\subset P_2}\eps_P^G k_{P,\of}^{\leq 0}(g).
$$
La proposition suivante va jouer un rôle auxiliaire analogue à celui de la proposition \ref{prop:cv}.

  \begin{proposition}\label{prop:cv-s}
L'intégrale
$$J_{P_1,\of}^{P_2,\sharp}(s)=\int_{P_1(F)\back G(\AAA)}F^{ P_1}(g) \sigma_1^2(H_{P_1}( g))k_{P_1,\of}^{P_2,\leq 0}(g) |\det(g)|^s\, dg
  $$
est absolument convergente pour $\Re(s)>0$ et définit une fonction holomorphe sur ce domaine. De plus, on a 
$$
\lim_{s\to 0}   s J_{P_1,\of}^{P_2,\sharp}(s)= \int_{P_1(F)\back G(\AAA)^1}F^{ P_1}(g) \sigma_1^2(H_{P_1}( g))k_{P_1,\of}^{P_2}(g) \, dg.
$$
\end{proposition}

\begin{preuve}
On suppose dans la preuve que $s>0$ est réel. Il est facile d'étendre ensuite le résultat à $s$ de partie réelle $>0$.
  Soit $(n_1',\ldots,n_l')$ le type de $P_2$. Soit $(n_{1},\ldots,n_{r_l})$ le type de $P_1$ avec pour $0\leq i\leq l-1$ et $r_0=0$
$$n_{r_i+1}+n_{r_i+2}+\ldots+ n_{r_{i+1}-1}=n'_{i+1}$$

Soit $g\in G(\AAA)$ et soit $g=mnk$ la décomposition d'Iwasawa de $g$ avec $m\in M_1(\AAA)$, $n\in N_1(\AAA)$ et $k\in K$. On écrit 
$$m=(m_1,\ldots,m_{r_l})
$$
avec $m_i\in GL(n_i,\AAA)$. Soit  $E'$ la fonction caractéristique des $g\in G(\AAA)$ tels que 
$$|\det(m_{r_i+1})|\leq 1
$$
pour $0\leq i<l$.

Soit $P$ un sous-groupe parabolique tel que $P_1\subset P\subset P_2$. Sur le lieu  des $g\in G(\AAA)$ tels que $F^{ P_1}(g) \sigma_1^2(H_{P_1}( g))=1$, les fonctions $E'$ et $E^P$ coïncident. Il s'ensuit que sur ce lieu on a 
$$
k_{P_1,\of}^{P_2,\leq 0}(g)=E^{'}(g) k_{P_1,\of}^{P_2}(g).
$$

Par conséquent,  en utilisant la décomposition $G(\AAA)=G(\AAA)^1 \times \RR_+^\times$, on a 
\begin{eqnarray*}
  J_{P_1,\of}^{P_2,\sharp}(s)&=& \int_{P_1(F)\back G(\AAA)} F^{ P_1}(g) \sigma_1^2(H_{P_1}( g))k_{P_1,\of}^{P_2}(g) E'(g)|\det(g)|^s  \, dg\\
&=& \int_{P_1(F)\back G(\AAA)^1} F^{ P_1}(g) \sigma_1^2(H_{P_1}( g))k_{P_1,\of}^{P_2}(g) c(g,s)  \, dg
\end{eqnarray*}

où l'on introduit
\begin{itemize}
\item $\displaystyle b(g)=\min_{0\leq i \leq l-1} |\det(m_{r_i+1})|^{-n/n_{r_i+1}}$ ;
\item  $\displaystyle  c(g,s)=  \int_0^{b(g)} a^{s-1} da= \frac{b(g)^s}s.$ (la mesure $da$ est la mesure de Lebesgue sur $\RR$).
\end{itemize}
Pour tout $g\in G(\AAA)^1$ tel que    $ \sigma_1^2(H_{P_1}( g))=1$, on a 
$$b(g)\leq 1\ ;$$
cette inégalité découle de l'inégalité ci-dessous, vraie sous cette même condition 
$$
1=\prod_{i=1}^{r_l} |\det(m_i)|\leq \prod_{i=0}^{l-1} |\det(m_{r_i+1})|^{n_{i+1}'/n_{r_i+1}}.
$$
On a donc pour  $g\in  G(\AAA)^1$ tel que    $ \sigma_1^2(H_{P_1}( g))=1$
$$|c(s,g)|\leq \frac1s.$$
La proposition résulte alors de la proposition \ref{prop:cv} et du théorème de convergence dominée.
\end{preuve}
\end{paragr}

\begin{paragr}[Une formule  limite.] --- Voici le principal résultat de cette section.
  
\begin{theoreme}\label{thm:holo}
  Pour tout $s\in \CC$ tel que $\Re(s)>0$, l'intégrale 
  \begin{equation}
    \label{eq:Jofsharp}
    J^\sharp_\of(f,s)=\int_{[G]} k^\sharp_\of(g) |\det(g)|^{s}\,dg
  \end{equation}
  converge absolument. En outre, on a 
$$
\lim_{s\to 0^+}  s J^\sharp(s)= J_\of(f).
$$
\end{theoreme}

\begin{remarque}
  Dans le théorème et dans tout l'article $\lim_{s\to 0^+}$ signifie que la limite est prise sur les $s\in \CC$ tels que $\Re(s)>0$.
\end{remarque}

\begin{preuve}
Comme dans la preuve du théorème \ref{thm:cv}, l'expression   $k^\sharp_\of(f, g)$ est égale à 
$$\sum_{P_1\subset P_2} \sum_{\delta\in P_1(F)\back G(F)}  F^{ P_1}(\delta g) \sigma_1^2(H_{P_1}( \delta g)) k_{P_1,\of}^{P_2,\leq 0}(f,\delta g) .$$
Le théorème résulte alors de la proposition \ref{prop:cv-s} et du théorème \ref{thm:cv}.
\end{preuve}

\end{paragr}

\section{Intégrale orbitale tronquée}

\begin{paragr}
On généralise légèrement la définition \eqref{eq:ksharp} : pour tout sous-groupe parabolique standard $Q$, on pose

\begin{equation}
  \label{eq:kQsh}
  k^{Q,\sharp}_\of(g)=\sum_{P\subset Q} \eps_P^Q \sum_{\delta \in P(F)\back G(F)}\hat{\tau}_P^Q(H_P(\delta g)) k_{P,\of}^{\leq 0}(\delta g),
\end{equation}
la somme portant sur les sous-groupes paraboliques standard $P$.
On a alors la formule d'inversion suivante.

\begin{lemme}\label{lem:inversion}On a 
$$
\sum_{Q}\sum_{\delta\in Q(F)\back G(F)} \tau_Q(H_Q(\delta g))k^{Q,\sharp}_\of(\delta g)=k_{G,\of}^{\leq 0}(g),
$$
où la somme porte sur les sous-groupes paraboliques standard $Q$.
\end{lemme}

\begin{preuve}
En utilisant la définition \ref{eq:kQsh}, on obtient

\begin{eqnarray*}
  \sum_{Q}\sum_{\delta\in Q(F)\back G(F)} \tau_Q(H_Q(\delta g))k^{Q,\sharp}_\of(\delta g)&=& \sum_{P\subset Q} \sum_{\delta \in P(F)\back G(F)}\eps_P^Q\hat{\tau}_P^Q(H_P(\delta g))\tau_Q(H_Q(\delta g))k_{P,\of}^{\leq 0}(\delta g) \\
&=& \sum_{P} \sum_{\delta \in P(F)\back G(F)}k_{P,\of}^{\leq 0}(\delta g) [\sum_{P\subset Q}\eps_P^Q\hat{\tau}_P^Q(H_P(\delta g))\tau_Q(H_Q(\delta g)) ]\\
&=& k_{G,\of}^{\leq 0}( g)
\end{eqnarray*}
car le crochet est nul sauf si $P=G$ auquel cas il vaut $1$ (c'est encore le lemme de Langlands \eqref{eq:langlands} utilisé dans la preuve de la proposition \ref{prop:ecriture}).
  \end{preuve}

\end{paragr}

\begin{paragr}[Une convergence auxiliaire.] 

  \begin{proposition}\label{prop:cv-s-Q}
     Pour tout sous-groupe parabolique standard $Q$ et toute fonction $f \in \Sc(\ggo(\AAA)) $, l'intégrale 
 $$\int_{Q(F)\back G(\AAA)} \tau_Q(H_Q(g)) k^{Q,\sharp}_\of(f,g)\, |\det(g)|^s dg$$
converge absolument pour $\Re(s)>0$. 
  \end{proposition}

   \begin{preuve} On prend $s>0$ réel.
   On introduit la variante suivante de \eqref{eq:sigma12}
     \begin{equation}
        \label{eq:sigma12Q}
    \sigma_{P_1}^{P_2,Q}(H)=\sum_{P_2\subset P\subset Q} \eps_{P_2}^{P}\tau_{P_1}^{P}(H) \hat{\tau}_{P}^Q(H) 
  \end{equation}
L'expression $k^{Q,\sharp}_\of(f,g)$ est alors égale à 
$$\sum_{P_1\subset P_2\subset Q} \sum_{\delta\in P_1(F)\back Q(F)}  F^{ P_1}(\delta g)  \sigma_{P_1}^{P_2,Q}(H_{P_1}( \delta g)) k_{P_1,\of}^{P_2,\leq 0}(f,\delta g) .$$
Il suffit donc de prouver la convergence pour $P_1\subset P_2\subset Q$ de 
$$\int_{P_1(F)\back G(\AAA)} \tau_Q(H_Q(g))  F^{ P_1}(g)  \sigma_{P_1}^{P_2,Q}(H_{P_1}( g))  |k_{P_1,\of}^{P_2,\leq 0}(f,g) |  \, |\det(g)|^s dg.
$$
Sous la condition  $F^{ P_1}(g)  \sigma_{P_1}^{P_2,Q}(H_{P_1}( g)) =1$, on  a, comme dans la preuve de la proposition \ref{prop:cv-s} dont on reprend les notations,
$$
k_{P_1,\of}^{P_2,\leq 0}(g)=E'(g) k_{P_1,\of}^{P_2}(g).
$$

Après décomposition d'Iwasawa $G(\AAA)=N_Q(\AAA)M_Q(\AAA)^1 A_Q^\infty K$, on est ramené à considérer l'intégrale
$$\int_{(M_Q\cap P_1)(F)\back M_Q(\AAA)^1 }  F^{ P_1}(m)  \sigma_{P_1}^{P_2,Q}(H_{P_1}( m))  \big(\int_K |k_{P_1,\of}^{P_2}(f,mk) |\,dk \big)   c_M(m,s)\, dm.
$$
où l'on introduit
$$c_M(m,s)= \int_{A_Q^\infty}  \tau_Q(H_Q(a)) E'(am)|\det(a)|^s \, da.$$ 
Soit $(n_1,\ldots,n_l)$ le type de $Q$ de sorte qu'on identifie $A_Q^\infty$ à $(\RR^\times_+)^l$. Soit $a=(a_1,\ldots,a_l)\in (\RR^\times_+)^l$. Par un raisonnement similaire à celui effectué  dans la preuve de la proposition \ref{prop:cv-s}, on voit que $E'(am)=1$ implique qu'on a $a_1 \leq 1$. Par conséquent, en utilisant la mesure de Lebesgue $da_i$, on majore  $c_M(m,s)$ par 
\begin{eqnarray*}
  \int_0^1 a_1^{s-1}\int_{0}^{a_1} a_2^{s-1} \ldots a_{l-1}^{s-1}\int_0^{a_{l-1}} a_l^{s-1}da_l \ldots da_1&=& \frac1s    \int_0^1 a_1^{s-1}\int_{0}^{a_1} a_2^{s-1} \ldots \int_0^{a_{l-2}}a_{l-1}^{2s-1}da_{l-1}\\
&=& \frac1{l!s^l}.
\end{eqnarray*}
On peut alors conclure parce que l'intégrale 
$$\int_{(M_Q\cap P_1)(F)\back M_Q(\AAA)^1 }  F^{ P_1}(m)  \sigma_{P_1}^{P_2,Q}(H_{P_1}( m))  \big(\int_K |k_{P_1,\of}^{P_2}(f,mk) |\,dk \big)\, dm
$$
est convergente.
   \end{preuve}

 \end{paragr}
 
 \begin{paragr}[Intégrale orbitale tronquée.] ---  Elle est définie par l'égalité \eqref{eq:IOtronq} ci-dessous.

   \begin{proposition}\label{prop:IOtronq}
L'intégrale
         \begin{equation}
           \label{eq:IOtronq}
           J_{\of}^{G,\leq 0}(f,s)= \int_{[G]}  k_{G,\of}^{\leq 0}( f,g) |\det(g)|^s \,dg
         \end{equation}
         converge absolument pour $\Re(s)>0$.
   \end{proposition}

   \begin{preuve}
     La proposition  est une conséquence directe de la formule d'inversion (lemme \ref{lem:inversion}) et de la proposition \ref{prop:cv-s-Q} ci-dessus.
   \end{preuve}
 \end{paragr}

\begin{paragr}\label{S:IOLevi}    Plus généralement, pour un sous-groupe de Levi $M$ de $G$,  $f\in \Sc(\mgo(\AAA))$ et $\of$ une orbite dans $\mgo(F)$, on définit 
$$k^{M,\leq 0}_\of(f,m)= E^{M}(m)\sum_{X\in \of(F)} f(m^{-1}Xm)$$
et 
$$J_{\of}^{M,\leq 0}(f,s)= \int_{[M]} k^{M,\leq 0}_\of(f,m)   |\det(m)|^s \,dm.
$$
De même, on montre que cette intégrale converge pour $\Re(s)>0$.
 \end{paragr}

\section{Contributions des orbites et intégrales orbitales}

 \begin{paragr}[Fonction presque invariante.] --- \label{S:pinv} Soit $f\in \Sc(\ggo(\AAA))$. Pour tout sous-groupe parabolique standard $P$ et $X\in \mgo_P(\AAA)$, on introduit le terme constant
$$f_P(X)=\int_{K} \int_{\ngo_P(\AAA)} f(k^{-1}(X+U)k)\,dUdk.$$
Cela définit un élément $f_P$ de $\Sc(\mgo_P(\AAA))$. Lorsque $M$ est un sous-groupe de Levi standard, on note $f_M=f_P$ où $P$ est le sous-groupe parabolique standard de facteur de Levi $M$.  Introduisons la définition suivante.

   \begin{definition}\label{def:inv}
     Une fonction $f\in \Sc(\ggo(\AAA))$ est dite \emph{presque invariante} si pour tous sous-groupes  paraboliques standard $P$ et $Q$ et tout $w\in W$ tel que  $M_Q=wM_Pw^{-1}$ on a
$$f_Q(X)= f_P(wXw^{-1})$$
pour tout $X\in \mgo_P(\AAA)$.
   \end{definition}

  \begin{remarque}\label{rq:inv}
À l'aide de l'injectivité de la transformation d'Abel sur les fonctions sphériques, il n'est pas difficile de voir qu'une fonction $K$-invariante à droite et à gauche est presque invariante au sens de la définition ci-dessus. Voici un exemple d'une telle fonction :  on prend $f=\otimes_{v} f_v$ telle que 
     \begin{itemize}
     \item en une place archimédienne $v$
$$f_v(X)=\varphi(\|X\|^2)$$
où $\varphi$ est une fonction  de Schwartz sur $\RR$ et $\|\cdot\|$ est la norme pour le produit hermitien standard sur $\ggo(F_v)$.
\item $f_v=\mathbf{1}_{\ggo(\oc_v)}$ en une place non-archimédienne $v$.
     \end{itemize}

   \end{remarque}
 \end{paragr}

 \begin{paragr} Soit $\of\in \oc$  et 
$$\lc(\of)$$
l'ensemble (fini) des couples $(M,\of')$ formé d'un  sous-groupe de Levi standard $M$ de $G$ et d'une orbite $\of'\in \oc^M$  telle que $I_M^G(\of')=\of$ où l'on note $I_M^G(\of')$ l'orbite induite comme dans  \cite{cuft} §2.9.  Deux éléments de $\lc(\of)$ sont équivalents s'ils sont conjugués par un élément du  groupe de Weyl $W$. Soit 
$$\lc(\of)/W$$
l'ensemble des classes d'équivalence qu'on identifie souvent  à un système de représentants.

Pour tout $(M,\of')\in \lc(\of)$ soit  $\stab_W(M,\of')$ son  stabilisateur dans $W$. Le groupe de Weyl $W^M$ relatif à $M$ est bien sûr un sous-groupe de  $\stab_W(M,\of')$ ; on note
\begin{equation}
  \label{eq:WMof}
  W(M,\of')=\stab_W(M,\of')/W^M
\end{equation}
le quotient. Pour tout sous-groupe de Levi de $G$ soit
\begin{equation}
  \label{eq:rM}
  r_M=(r-1)!
\end{equation}
où $r$ est  le nombre de \og blocs linéaires\fg{} de $M$.

 \begin{proposition}\label{prop:explicit} Soit $f\in \Sc(\ggo(\AAA))$ presque invariante au sens de la définition \ref{def:inv}. 
On  a 
$$J_{\of}^{G,\sharp}(f,s)=\sum_{(M,\of') \in \lc(\of)/W}\eps_M^G \frac{r_M}{|W(M,\of')|} J^{M,\leq 0}_{\of'}(f_M,s)
$$
\end{proposition}

Avant de donner la preuve de cette proposition, notons le corollaire suivant qui repose sur le théorème \ref{thm:holo}.

\begin{corollaire}\label{cor:explicit}
  Sous les hypothèses de la proposition 
 $$J_\of(f)=\lim_{s\to 0^+}  s  \sum_{(M,\of') \in \lc(\of)/W}\eps_M^G \frac{r_M}{|W(M,\of')|}  J^{M,\leq 0}_{\of'}(f_M,s).
$$
\end{corollaire}

\end{paragr}

\begin{paragr}[Démonstration de la proposition \ref{prop:explicit}.] ---  On va faire une série de manipulations formelles. Il est aisé de les justifier  pour $\Re(s)>0$, par exemple en majorant $f$ par une fonction positive et en prenant $s$ réel. En utilisant successivement la formule \eqref{eq:kQsh} pour $Q=G$ (qui n'est autre que la définition \eqref{eq:ksharp}),  la décomposition d'Iwasawa combinée à un changement de variables, on a 
\begin{eqnarray*}
  J_{\of}^{G,\sharp}(f,s)&=&\sum_{P} \eps_P^G \int_{P(F)\back G(\AAA)} \hat{\tau}_P(H_P(g)) k_{P,\of}^{\leq 0}(f,g) \,|\det(g)|^sdg\\
&=&\sum_{P} \eps_P^G \sum_{\{\of'\in \oc^{M_P} \mid I_P^G(\of')=\of \}}\int_{[M_P]} \hat{\tau}_P(H_P(m)) k_{\of'}^{M_P,\leq 0}(f_P,m) \,|\det(m)|^sdm,
\end{eqnarray*}
où  
$$k_{\of'}^{M_P,\leq 0}(f_P,m)=E^{M_P}(m) \sum_{X\in \of'(F)} f_P(m^{-1}Xm).$$
Soit $(M,\of')\in \lc(\of)$. Soit
$$W_M^+=\{w\in W \mid w\Delta_0^M\subset \Delta_0\}.
$$
L'application $w\mapsto (wMw^{-1},w\cdot\of')$ fournit une surjection de $W_M^+$ sur la $W$-orbite de $(M,\of')$. Le cardinal d'une fibre est exactement $|W(M,\of')|$. Pour tout $w\in W_M^+$ posons 
 $P_w=Mw^{-1}Bw$. C'est un sous-groupe parabolique de facteur de Levi $M$ c'est-à-dire un élément de $\pc(M)$. L'application
$$w\mapsto P_w
$$
induit une bijection de $W^+_M$ sur $\pc(M)$. 

Soit $P$ et $P_1$ des sous-groupes paraboliques standard de facteurs de Levi respectifs $M$ et $M_1$. On complète ces données en des éléments  $(M,\of')$ et $(M_1,\of_1')$ de $\lc(\of)$. Supposons que $w\in W$ vérifie  $(M_1,\of_1')=w\cdot (M,\of')$. En utilisant le fait que $f$ est presque invariante, on obtient par un changement de variables qu'on a
\begin{eqnarray*}
  \int_{[M_1]} \hat{\tau}_{P_1}(H_{M_1}(m)) k_{\of'_1}^{M_1,\leq 0}(f_{M_1},m) \,|\det(m)|^sdm\\
= \int_{[M]} \hat{\tau}_{P_w}(H_{M}(m)) k_{\of'}^{M,\leq 0}(f_{M},m) \,|\det(m)|^sdm
\end{eqnarray*}

On en déduit qu'on a 
$$
  J_{\of}^{G,\sharp}(f,s)=\sum_{(M,\of') \in \lc(\of)/W }\eps_M^G |W(M,\of')|^{-1}\int_{[M]}  \big(\sum_{w\in W_M^+} \hat{\tau}_{P_w}(H_{M}(m)) \big)k_{\of'}^{M,\leq 0}(f_{M},m) \,|\det(m)|^s dm
$$
où 
$$\sum_{w\in W_M^+} \hat{\tau}_{P_w}(H_{M}(m))=\sum_{P\in \pc(M)} \hat{\tau}_P(H_{M}(m)).$$  D'après le lemme \ref{lem:sumtau} ci-dessous, on peut remplacer dans l'intégrale cette expression par $r_M$. Le résultat s'ensuit.

\begin{lemme}
  \label{lem:sumtau}
Pour tout $H\in a_M$ en dehors d'un ensemble de mesure nulle, on  a 
$$
 \sum_{P\in \pc(M)}  \hat{\tau}_{P}(H)=r_M.
$$
\end{lemme}

\begin{preuve}
On raisonne par récurrence sur la dimension de $a_M$. D'après le lemme de Langlands, pour tout sous-groupe parabolique $P$, on 
$$\sum_{P\subset Q\subset G }  \eps_Q^G \tau_P^Q(H) \hat{\tau}_Q^G(H)=0.$$
En sommant sur $P\in \pc(M)$ et en intervertissant les sommes sur $P$ et $Q$ on aboutit à 
$$\sum_{Q\in \fc(M)}  \eps_Q^G \hat{\tau}_Q^G(H) [ \sum_{P\in \pc^Q(M)} \tau_P^Q(H)] =0.$$
Le facteur entre crochets vaut $1$ presque partout. L'égalité devient
$$\sum_{L\in \lc(M)} \eps_L^G  [\sum_{Q\in \pc(L)} \hat{\tau}_Q^G(H)].$$
Par hypothèse de récurrence, le terme entre crochets  pour $L\not=M$ est égal à $(k-1)!$ où $k$ est le nombre de blocs de $L$. Du coup l'égalité se réécrit
\begin{equation}
  \label{eq:stirling}
 (-1)^{N-1} \sum_{P\in \pc(M)} \hat{\tau}_P^G(H)  +\sum_{k=1}^{N-1}  (-1)^{k-1}  (k-1)! |\lc_k(M)|=0
\end{equation}
où $N$ est le nombre de blocs de $M$ et où $|\lc_k(M)|$ est le cardinal de l'ensemble des éléments de $\lc(M)$ qui sont formés de $k$ blocs. On peut encore interpréter $|\lc_k(M)|$ comme le nombre de partitions en $k$ ensembles d'un ensemble à $N$ éléments. On a alors la relation bien connue
$$
\sum_{k=1}^N |\lc_k(M)| X(X-1)\cdots(X-k+1) = X^N
$$
qui donne après division par $X$ et évaluation en $X=0$
\begin{equation}
  \label{eq:stirling2}
\sum_{k=1}^N |\lc_k(M)| (-1)^{k-1} (k-1)!=0.
\end{equation}
Comme $|\lc_N(M)|=1$, on obtient le résultat cherché en comparant \eqref{eq:stirling} et \eqref{eq:stirling2}.
\end{preuve}
\end{paragr}

\section{Intégrales tronquées de séries d'Eisenstein}\label{sec:Eis}

\begin{paragr} Dans toute la suite,  $\of\in \oc$ est une orbite nilpotente.
  
\end{paragr}

\begin{paragr}[Fonction $\theta_\of(s)$.] ---  \label{S:theta} Soit $d_i$ la multiplicité du bloc de taille $i$ dans la décomposition de Jordan de $\of$. Soit $d=\sum_{i=0}^\infty d_i$ et 
$$r=\max\{j\mid d_j\not=0\}.
$$

 Soit $P$ le sous-groupe parabolique standard de $GL(d)$ de type $(d_r,d_{r-1},\ldots, d_1)$. On écrit un élément $m\in M_P$ sous la forme $((m_r,m_{r-1},\ldots,m_1)$ avec $m_i\in GL(d_i)$. Soit $\la_P\in a_P^*$ tel que  pour tout $m\in M_P(\AAA)$ on ait
$$\exp(\bg \la_P,H_P(m))= \prod_{i=1}^r |\det(m_i)|^{i}.
$$
On définit alors
$$\ec_\of(s,g)=\sum_{\delta\in P(F)\back G(F)} \exp(\bg 2\rho_P+s\la_P,H_P(\delta g))
$$
où $\rho_P$ est la demi-somme des racines dans $N_P$. La théorie élémentaire des séries d'Eisenstein assure que cette série converge pour $\Re(s)>0$ et définit sur cet ouvert une fonction holomorphe de la variable $s$.

\begin{proposition}\label{prop:theta}
L'intégrale 
$$\theta_\of^G(s)=\int_{[G]} E(g) \ec_\of(s,g) \, dg$$
converge absolument pour $\Re(s)>0$ et définit sur cet ouvert une fonction holomorphe.
\end{proposition}

\begin{preuve} Pour obtenir les majorations nécessaires pour obtenir ce résultat on suppose que $s>0$ est réel. 
  En utilisant \eqref{eq:EQ}, on est conduit à étudier  l'intégrale
$$\int_{Q(F)\back G(\AAA)} F^Q(g) \tau_Q(H_Q(g)) \chi_Q^G(H_Q(g)) \ec_\of(s,g) \,dg.$$
pour $Q$ un sous-groupe parabolique standard. Soit $Q=MN$ la décomposition de Levi standard. Par décomposition d'Iwasawa, l'intégrale ci-dessus s'écrit
$$
\int_{[M]} \exp(-\bg2\rho_Q,H_Q(m)\bd) F^Q(m)  \tau_Q(H_Q(m))  \chi_Q^G(H_Q(m)) (\ec_\of)_Q(s,m) \,dm,
$$
où $(\ec_\of)_Q(s,m) $ désigne le terme constant de la série d'Eisenstein $\ec_\of(s,m)$ le long de $Q$. En utilisant le calcul classique de ce terme constant, on est ramené à considérer l'intégrale pour $w\in W$ tel que $R=M \cap w^{-1}Pw$ est un sous-groupe parabolique standard de $M$
$$
\int_{[M]} \exp(-\bg2\rho_Q,H_Q(m)\bd) F^Q(m)  \tau_Q(H_Q(m))  \chi_Q^G(H_Q(m)) \ec_w(s,m) \,dm,
$$
où l'on introduit la série convergente
$$\ec_w(s,m) =\sum_{\delta \in R(F)\back M(F)} \exp(\bg 2\rho_R^M  +\sigma_w   +sw^{-1}\la_P),H_R(\delta m)\bd)$$
avec $\sigma_w$ est la somme des racines dans $N_Q\cap w^{-1}N_Pw$.
Cette intégrale est donc le produit des deux intégrales
\begin{equation}
  \label{eq:integ1}
  \int_{M(F)\back M(\AAA)^1} F^Q(m)   \ec_w(s,m) \,dm
\end{equation}
et
\begin{equation}
  \label{eq:integ2}
  \int_{A_M^\infty}   \exp(\bg -2\rho_Q  +2\rho_R^M  +\sigma_w   +sw^{-1}\la_P  ,H_M(a)\bd)  \tau_Q(H_Q(a))  \chi_Q^G(H_Q(a)) \,dm.
\end{equation}
L'intégrale \eqref{eq:integ1} converge car la restriction de $F^Q$ à $M(F)\back M(\AAA)^1$ est la fonction caractéristique d'un compact. L'intégrale \eqref{eq:integ2} est élémentaire. On peut observer que la projection sur $a_M^*$ de   $ -2\rho_Q  +2\rho_R^M  +\sigma_w$ est une somme à coefficients négatifs  d'éléments de $\Delta_Q$. Par ailleurs, on vérifie aisément que la projection  sur $a_M^*$ de  $w^{-1}\cdot\la_P$ est une combinaison linéaire à coefficients négatifs de l'unique élément de $\varpi_Q$  et  d'éléments de $\Delta_Q$. La convergence s'ensuit.
\end{preuve}

\begin{remarque}\label{rq:calcul-theta}
  La démonstration de la proposition \ref{prop:theta} donne en même maintenant un moyen d'expliciter les fonctions $\theta_\of^G(s)$. Il suffit essentiellement d'évaluer l'intégrale \eqref{eq:integ1} ci-dessus. En fait, le cas essentiel, auquel on se ramène par un calcul de résidu, est l'intégrale
$$\int_{G(F)\back G(\AAA)^1} F^G(g) \ec(g,\la)\,dg$$
où 
$$\ec(g,\la)=\sum_{\delta \in B(F)\back G(F)} \exp(\bg 2\rho_B +\la,H_B(\delta g)\bd).$$
En utilisant l'opérateur de troncature $\Lambda$ d'Arthur( cf. \cite{ar2} en le paramètre $T=0$), on est ramené à évaluer l'intégrale
 $$\int_{G(F)\back G(\AAA)^1} \Lambda \ec(g,\la)\,dg.$$
La fonction  $\Lambda \ec(g,\la)$ se calcule très simplement à l'aide des opérateurs  d'entrelacement (cf. \cite{Ar-intro} lemme 15.2). Ceux-ci s'explicitent très bien  dans cette situation. Le calcul de l'intégrale s'en déduit (pour une méthode légèrement différente, cf. \cite{KW} section 3).  
\end{remarque}
\end{paragr}

\begin{paragr}[Généralisation.] ---  Plus généralement on définit une fonction $\theta_{\of}^M(s)$ pour $M=GL(n_1)\times \ldots \times GL(n_r)$ un sous-groupe de Levi de $G$ et $\of'$ une orbite nilpotente dans $\mgo$ (c'est-à-dire une collection d'orbites nilpotentes $\of_i$ pour $GL(n_i)$ et $1\leq i\leq r$)  de la façon suivante.
$$\theta_{\of'}^M(s)=\prod_{i=1}^r\theta^{GL(n_i)}_{\of_i}(s).$$

\end{paragr}

\section{Intégrales zêta associées à une orbite nilpotente}\label{sec:zeta}

\begin{paragr}
  On continue avec les notations de la section  \ref{sec:Eis}.
\end{paragr}

\begin{paragr}  \label{S:IO}   On associe à l'orbite $\of$ le sous-groupe parabolique standard $R$ de $G$ de type 
$$(d_r,d_{r_1}, \ldots,d_1, d_r,\ldots,d_2,\ldots, d_r,d_{r-1},d_r).$$
Soit $R=M_R N_R$ sa décomposition de Levi.  Soit $X\in \of$ l'élément défini matriciellement par
 $$
\begin{pmatrix} 0_{d_1+\ldots +d_r} & \underset{0_{}}{I_{d_2+\ldots+d_r  }} & 0 & 0&  0 \\   &0_{d_2+\ldots+d_r}  & 0 & 0& 0\\ & & \ddots & \underset{0}{I_{d_{r-1}+d_r}} & 0\\ & & & 0_{d_{r-1}+d_r} & \underset{0}{I_{d_{r}}} \\ & & & & 0_{d_r}\end{pmatrix}
$$
La $N_R$-orbite de $X$ est de la forme $X+\ugo_{\of}$ où  $\ugo_{\of} \subset \ggo$ est un sous-espace de description simple qu'on ne rappellera pas ici (cf. proposition 4.5.1 de \cite{scuft}).

Soit $f\in \Cc(\ggo(\AAA))$ et 
$$f^K(X)=\int_K f(k^{-1}Xk).$$
Soit $L= \prod_{j=1}^r \prod_{i=1}^{j-1} GL(d_j)$. Un élément $A$ de $L$ va s'écrire  $(A_{i,j})_{ 1\leq i <j\leq r}$ avec $A_{i,j}\in GL(d_j)$. On pose pour $A\in L(\AAA)$
$$f^K_{\of}(A)=\int_{\ugo_{\of}(\AAA)} f^K( U+   \begin{pmatrix} 0_{d_1+\ldots +d_r} & \underset{0_{}}{\Delta_1(A)} & 0 & 0&  0 \\   &0_{d_2+\ldots+d_r}  & 0 & 0& 0\\ & & \ddots & \underset{0}{\Delta_{r-2}(A)} & 0\\ & & & 0_{d_{r-1}+d_r} & \underset{0}{\Delta_{r-1}(A)} \\ & & & & 0_{d_r}\end{pmatrix}) dU$$
où pour $1\leq i \leq r-1$
 $$\Delta_i(A)= \begin{pmatrix} A_{i,r} & * & * &* \\  & A_{i,r-1}  & * &* \\ & & \ddots & * \\ & & & A_{i,i+1}
  \end{pmatrix}
  $$
La mesure de Haar sur $\ugo_\of(\AAA)$ est normalisée de sorte que $\ugo_\of(F)$ soit de covolume $1$.
\end{paragr}

\begin{paragr}[Fonctions zêta.] ---\label{S:zeta}
  On introduit alors la fonction zêta
$$
Z_\of^G(f,s)= \int_{L(\AAA)}  f^K_{\of}(A)  \delta(A,s) dA 
$$
où
$$\delta(A,s)= \prod_{j=1}^r \prod_{i=1}^{j-1} |\det(A_{i,j})|^{d_i+\ldots+d_j+(j-i)s  }.
$$
La normalisation de la mesure de Haar sur $L(\AAA)$ suit les prescriptions du §\ref{S:notations}.

\begin{proposition}\label{prop:zeta}
  L'intégrale $Z_{\of}^G(s,f)$ converge pour $\Re(s)>0$ et elle définit une fonction holomorphe sur cet ouvert.
\end{proposition}

\begin{preuve}
  Elle est élémentaire (cf. \cite{GJ} ou \cite{scuft} §10.2).
\end{preuve}
\end{paragr}

\begin{paragr}[Généralisation aux sous-groupes de Levi.] --- Plus généralement, pour $M$  un sous-groupe de Levi de $G$, tout $f\in \Sc(\mgo(\AAA))$ et toute $M$-orbite nilpotente $\of$ dans $\mgo$, on définit une intégrale zêta $Z^M_\of(f)$. Si $f$ est un produit sur les blocs linéaires alors l'intégrale zêta est un produit des fonctions zêta associés aux blocs.  
  
\end{paragr}

\section{Calcul d'intégrales orbitales nilpotentes tronquées}

\begin{paragr} Le résultat suivant explicite  l'intégrale orbitale tronquée $J^{G,\leq 0}_{\of}(f,s)$ en termes d'une intégrale d'une série d'Eisenstein introduite à la section \ref{sec:Eis} et de la distribution zêta introduite à la section \ref{sec:zeta}.

  \begin{theoreme}\label{thm:J0-J}
    Soit $f\in \Cc(\ggo(\AAA))$.  Si le support de $f$ est \og suffisamment petit\fg{}, alors on a pour tout $s\in\CC$ tel que $\Re(s)>0$
$$J^{G,\leq 0}_{\of}(f,s)= Z_{\of}^G(f,s)\cdot \theta_{\of}^G(s).
$$
  \end{theoreme}

 \begin{remarque}
    Soit $f\in \Cc(\ggo(\AAA))$ une fonction à support compact. On fixe une place $v$ et $t\in F_v$. Alors le support de
$$f(t^{-1}\cdot)
$$
est suffisamment petit au sens du théorème \ref{thm:J0-J} dès que $|t|_v$ est assez petit. Les intégrales d'Arthur $J_\of(f)$ vérifient des propriétés d'homogénéité (cf. \cite{ar_unipvar}). En principe, cette hypothèse de petitesse du support n'est pas trop restrictive.
  \end{remarque}

\begin{preuve}
Soit $Q$ le sous-groupe parabolique standard de type $(d_1+\ldots+d_r,\ldots,d_{r-1}+d_r,d_r)$. Soit $X\in \of$ l'élément défini au §\ref{S:IO}. Soit $G_X$ le centralisateur de $X$ dans $G$. On a $G_X\subset Q$. L'intégrale s'écrit par décomposition d'Iwasawa 

$$J^{G,\leq 0}_{\of}(f,s)=\int_{G_X(F)\back Q(\AAA)} \int_K f(k^{-1} q^{-1}Xqk) \, dk E^G(q) |\det(q)|^s \,dq.$$
Soit $Q=MN$ la décomposition de Levi standard de $Q$. Soit $q=mn$ avec $m\in M(\AAA)$ et $n\in N(\AAA)$. On écrit encore $m=(m_1,\ldots,m_r)$ avec $m_i\in GL(d_i+\ldots+d_r,\AAA)$. On veut montrer que sous la condition 
\begin{equation}
  \label{eq:intnonnulle}
   \int_K f(k^{-1} q^{-1}Xqk) \, dk \not=0
 \end{equation}
 et le fait que le support de $f$ est assez petit qu'on a l'équivalence
$$E^G(q)=1 \Leftrightarrow E^{GL(d_1+\ldots+d_r)}(m_1).$$
Dans le contexte des fibrés de Higgs, cette propriété clef  a été observée par Schiffmann (cf. \cite{Schif}).
Le sens $\Rightarrow$ est toujours vrai (cf. corollaire \ref{cor:encad}). Supposons $E^{GL(d_1+\ldots+d_r)}(m_1)=1$. Sous la condition \eqref{eq:intnonnulle}, l'élément
\begin{equation}
  \label{eq:element}
m_1^{-1} \left(\underset{0_{}}{I_{d_2+\ldots+d_r  }} \right)m_2
\end{equation}
est astreint à rester dans un compact qui dépend du support de $f$. Soit $P$ un sous-groupe parabolique maximal de $GL(d_2+\ldots+d_r)$ et $\delta \in GL(d_2+\ldots+d_r,F)$. On regarde  $GL(d_2+\ldots+d_r)$ naturellement comme un sous-groupe de  $GL(d_1+\ldots+d_r)$ de sorte que l'élément s'écrit encore 
\begin{equation}
  \label{eq:element2}
(\delta m_1)^{-1} \left(\underset{0_{}}{I_{d_2+\ldots+d_r  }} \right)(\delta m_2). 
\end{equation}
Soit $(l,d_2+\ldots+d_r-l)$ le type de $P$. Soit $P'$ le sous-groupe parabolique standard de $GL(d_1+\ldots+d_r)$ de type $(l,d_1+\ldots+d_r-l)$. Soit $\delta m_2=pk$ et $\delta m_1=p'k'$ les décomposition d'Iwasawa relatives à $P$ et $P'$. On écrit
$$p=
\begin{pmatrix}
  a & * \\ 0& *
\end{pmatrix}
 \ \ \ \
p'=
\begin{pmatrix}
  b& * \\ 0& *
\end{pmatrix}
$$
avec $a$ et $b$ dans $GL(l,\AAA)$. En multipliant \eqref{eq:element2} à gauche par $k'$ et à droite par $k^{-1}$, on voit que l'élément
$$
\begin{pmatrix}
  b^{-1}a & * \\ 0& *
\end{pmatrix}
$$
est astreint à rester dans un compact qui ne dépend que du support de $f$. Si  celui-ci est assez petit alors on a nécessairement (pour la valeur absolue adélique)
$$|\det(b^{-1}a)|\leq 1.$$
 Or $|\det(b)|\leq 1$ (cf. corollaire \eqref{cor:car} qu'on utilise pour $m_1$ et la condition $E^{GL(d_1+\ldots+d_r)}(m_1)=1$). Il s'ensuit qu'on a $|\det(a)|\leq 1$ c'est-à-dire $\bg\varpi_P,H_P(\delta m_2)\bd\leq 0$. Comme on a cette inégalité pour tout $P$ et $\delta$, une nouvelle application du corollaire \ref{cor:car} entraîne qu'on a $E^{GL(d_2+\ldots+d_r)}(m_2)=1$. Le même raisonnement donne par récurrence 
$$E^{GL(d_i+\ldots+d_r)}(m_i)=1$$
pour $1\leq i \leq r$. Mais ces conditions  impliquent $E^G(q)=1$ (cf. corollaire \ref{cor:encad}).

On peut donc écrire 

$$J^{G,\leq 0}_{\of}(f,s)=\int_{G_X(F)\back Q(\AAA)} \int_K f(k^{-1} q^{-1}Xqk) \, dk E^H(m_1) |\det(q)|^s \,dq$$
où $H=GL(d_1+\ldots+d_r)$.
Soit $R\subset Q$ le sous-groupe parabolique de $G$ défini au §\ref{S:IO}. On écrit un élément $m$ de $M_R$ sous la forme $(m_r^1,m_{r-1}^1, \ldots, m_1^1, m_r^2,\ldots, m_2^2,m_3^r, \ldots, m_r^r)$. Avec cette écriture, on  a
$$m_1=(m_r^1,m_{r-1}^1, \ldots, m_1^1).$$

Soit $M_X\subset M_R$  le facteur de Levi de $G_X$ qui s'identifie à $GL(d_r)\times \ldots \times GL(d_1)$, chaque facteur $GL(d_j)$ étant plongé diagonalement.

Soit $N^H_R=H\cap N_R$ (on voit naturellement $H$ comme un sous-groupe de $G$). La décomposition d'Iwasawa et les lemmes 7.3.3 et  7.3.4 de \cite{scuft} impliquent qu'on a 
$$J^{G,\leq 0}_{\of}(f,s)=\int_{ M_X(F)\back M_R(\AAA)}   f_\of^K(  ((m_j^i)^{-1}m_j^{i+1})  )    \int_{[N_R^H]}E^{H}( nm_1)\,dn    \,  \delta(m)  |\det(   m)|^s  \, dm$$
où l'on pose
$$\delta(m)=  \prod_{j=1}^r \prod_{i=1}^{j-1}   |\det ((m_j^i)^{-1}m_j^{i+1}) |^{d_i+\ldots d_j}  .$$
On effectue alors le changement de variables $A_{i,j}=(m_j^i)^{-1}m_j^{i+1}$ et l'on pose $M^H_R$ ce qui donne
$$J^{G,\leq 0}_{\of}(f,s)= Z_{\of}(f,s) \cdot \int_{ [M_R^H]  } \int_{[N_R^H]} E^H(nm_1) \prod_{j=1}^r |\det(m_j^1)|^{js} \,dndm_1.
$$
On reconnaît sans mal dans le second facteur du membre de droite  la fonction $\theta^G_\of(s)$. Cela conclut.

  \end{preuve}
\end{paragr}

\section{L'énoncé final}

\begin{paragr} En utilisant le corollaire \ref{cor:explicit} et la généralisation évidente du théorème \ref{thm:J0-J} aux sous-groupes de Levi de $G$, on aboutit à l'énoncé suivant.
  
  \begin{theoreme} \label{thm:residu}
Soit $f\in \Cc(\ggo(\AAA))$ presque invariante au sens de la définition \ref{def:inv}. Lorsque le support de $f$ est assez petit, on a 
$$J_\of(f)=\lim_{s\to 0^+}s  \sum_{(M,\of') \in \lc(\of)/W}\eps_M^G \frac{r_M}{|W(M,\of')|}Z_{\of'}^M(f_M,s)\cdot \theta_{\of'}^M(s).
$$
  \end{theoreme}
\end{paragr}

\section{Une conjecture}

\begin{paragr}
 Désormais, on prend $F=\QQ$. On aurait pu rester avec un corps de nombres général mais la restriction à $\QQ$ simplifie légèrement les formules. Soit pour $s\in \CC$
$$\xi(s)=\pi^{-s/2} \Gamma(s/2) \zeta_\QQ(s)$$
le produit de la fonction zêta $\zeta_\QQ$ de Riemann avec le facteur archimédien usuel, $\Gamma$ étant la fonction d'Euler.
\end{paragr}

\begin{paragr} Soit $f\in \Sc(\ggo(\AAA))$ définie par $f=f_\infty\otimes (\otimes_p f_p)$ où 
 \begin{itemize}
  \item $f_\infty\in \Sc(\ggo(\RR))$ est donnée par $f_\infty(X)=\exp(-\pi \trace(\,^tXX))$, 
  \item pour $p$ premier, $f_p\in \Cc(\ggo(\Qp))$ est la fonction caractéristique de $\ggo(\ZZ_p)$.
  \end{itemize}
  Le but de cette section est de donner une formule conjecturale pour la distribution $J_\of(f)$ pour une orbite nilpotente $\of$. Pour cela, il nous faut normaliser la mesure sur $G(\AAA)^1$. Avec les conventions sur les choix de mesures du §\eqref{S:notations}, il suffit de préciser la mesure de Haar sur le tore standard $T_0(\AAA)$. Celui-ci est naturellement identifié à $(\AAA^1 \times \RR_+^\times)^n$ où $\AAA^1$ désigne le groupe multiplicatif des idèles de module $1$. On munit  $\RR_+^\times$ de la mesure $dt/t$ et $\AAA^1$ de la mesure qui donne le volume $1$ au quotient $\QQ^\times\back \AAA^1$. Avec ces choix, un calcul classique montre qu'on a 
  \begin{equation}
    \label{eq:volume}
    \vol(G(\QQ)\back G(\AAA)^1)=\xi^*(1)\xi(2)\cdots\xi(n)
  \end{equation}
  où $\xi^*(1)$ est le résidu en $s=1$ de $\xi$ (ici égal à $1$ mais c'est plus suggestif de le laisser sous cette forme).

\end{paragr}

\begin{paragr}[Fonctions zêta associés à une orbite nilpotente.] --- Soit $\of$ une orbite nilpotente dans $\ggo$. On reprend les notations du §\ref{S:theta}.  On lui associe le diagramme de Young $Y_\of$ suivant : il a d'abord $d_r$ colonnes chacune ayant $r$ cases, puis $d_{r-1}$ colonnes  de $r-1$ cases etc. La première ligne a alors $d$ cases.  Dans cette écriture, l'orbite nulle a pour diagramme une ligne de $n$ cases et l'orbite régulière une colonne de $n$ cases. On pose alors pour $s\in \CC$
  $$Z_\of^G(s)=\prod_{y\in Y_\of} \xi(1+b_y+c_ys)$$
où l'on introduit pour $y\in Y_\of$ les fonctions 
\begin{itemize}
\item bras $b_y$ qui compte le nombre de cases strictement à droite de $y$ ;
\item jambe $j_y$ qui compte le nombre de cases strictement en dessous  de $y$;
\item crochet $c_y=1+b_y+j_y$.
\end{itemize}
Ainsi, on obtient $\xi(1+s)\xi(2+2s)\cdots \xi(n+ns)$ pour l'orbite nulle, $\xi(1+s)\xi(1+2s)\cdots \xi(1+ns)$ pour l'orbite régulière et 
$\xi(1+s)^2 \xi(2+3s)$ pour l'orbite sous-régulière de $GL(3)$.
\end{paragr}

\begin{paragr}
  Soit $(M,\of')\in \lc(\of)$. On écrit $M=GL(n_1)\times \ldots \times GL(n_r)$ et $\of'=(\of_1',\ldots,\of_r')$. On généralise la définition précédente en posant
$$Z_{\of'}^{M}(s)=\prod_{i=1}^r Z^{GL(n_i)}_{\of'_i}(s).$$ 
Bien sûr cette expression ne dépend que de l'orbite sous $W$ de $(M,\of')$. On pose alors 
\begin{equation}
  \label{eq:Hof}
  H_\of(s)=\sum_{(M,\of') \in \lc(\of)/W}\eps_M^G \frac{r_M}{|W(M,\of')|}Z_{\of'}^M(s).
\end{equation}

On peut alors énoncer la conjecture.

\begin{conjecture}
  La fonction méromorphe $H_\of(s)$ a  un pôle simple en $s=0$ de résidu $J_\of(f)$. 
\end{conjecture}

Cette conjecture est un pendant pour le corps $\QQ$ d'une conjecture formulée sur les corps de fonctions dans \cite{scfh}. Cette dernière est elle-même un raffinement de conjectures de Hausel-Rodriguez-Villegas et Mozgovoy (cf. \cite{HRV} et \cite{Moz}). Elle a une certaine ressemblance formelle avec le théorème \ref{thm:residu}. Elle n'est évidente que lorsque $\of=(0)$ est l'orbite nulle : dans ce cas, on a simplement
$$H_\of(s)= \xi(1+s)\xi(2+2s)\cdots \xi(n+ns).$$
Le résultat est évident car d'une part on a   $J_{(0)}(f)= \vol(G(F)\back G(\AAA)^1) f(0)$ et d'autre part le volume est donné explicitement par \eqref{eq:volume}. En utilisant les résultats de \cite{scuft} et \cite{cuft}, on peut vérifier cette conjecture pour $n\leq 3$ (les calculs sont peu ou prou ceux de \cite{scfh} §8.5 et \emph{infra}). Par exemple, pour $\of$  l'orbite sous-régulière de $GL(3)$,  on a 
$$H_\of(s)=\xi(1+s)^2 \big(\xi(2+3s)- \xi(2+2s)\big)$$
qui a clairement un pôle simple en $s=0$ de résidu $\xi^*(1)^2\xi'(2)$.

Pour faire un peu plus le lien avec les constructions de \cite{scfh}, introduisons les sommes formelles d'orbites nilpotentes 
$$\sum_{\of} a_\of \of$$
où $\of$ parcourt l'ensemble des orbites nilpotentes dans $\mathfrak{gl}(n)$ quand $n$ décrit $\NN$ et $a_{\of}$ est une fonction méromorphe sur $\CC$. On a une structure additive évidente et une structure multiplicative induite par l'induction d'orbite : ainsi $\of_1\cdot \of_2$ est l'induite à $\mathfrak{gl}(n_1+n_2)$ de l'orbite $\of_1\times \of_2$ de $\mathfrak{gl}(n_1)\times \mathfrak{gl}(n_2)$. On a alors l'énoncé combinatoire suivant qu'on laisse en exercice au lecteur. On pose $Z_{(0)}=1$ et $H_{(0)}=0$ si $(0)$ est l'orbite de $\mathfrak{gl}(0)$.

\begin{proposition}
  On a l'identité
$$\log(\sum_{\of} Z_\of \of )=\sum_{\of} H_\of \of.$$
\end{proposition}

\end{paragr}

\bibliography{bibiog}
\bibliographystyle{plain}

\begin{flushleft}
Pierre-Henri Chaudouard \\
Université Paris Diderot (Paris 7) et Institut Universitaire de France\\
 Institut de Mathématiques de Jussieu-Paris Rive Gauche \\
 UMR 7586 \\
 Bâtiment Sophie Germain \\
 Case 7012 \\
 F-75205 PARIS Cedex 13 \\
 France
\medskip

Adresse électronique :\\
Pierre-Henri.Chaudouard@imj-prg.fr \\
\end{flushleft}

\end{document}